
\magnification1200
\input amssym.def 
\input amssym.tex 
\def\SetAuthorHead#1{}
\def\SetTitleHead#1{}
\def\NoindentAfter{\everypar={\setbox0=\lastbox\everypar={}}}
\def\H#1\par#2\par{{\baselineskip=15pt\parindent=0pt\parskip=0pt
 \leftskip= 0pt plus.2\hsize\rightskip=0pt plus.2\hsize
 \bf#1\unskip\break\vskip 4pt\rm#2\unskip\break\hrule
 \vskip40pt plus4pt minus4pt}\NoindentAfter}
\def\HH#1\par{{\bigbreak\noindent\bf#1\medbreak}\NoindentAfter}
\def\HHH#1\par{{\bigbreak\noindent\bf#1\unskip.\kern.4em}}
\def\th#1\par{\medbreak\noindent{\bf#1\unskip.\kern.4em}\it}
\def\endth{\medbreak\rm}
\def\pf#1\par{\medbreak\noindent{\it#1\unskip.\kern.4em}}
\def\df#1\par{\medbreak\noindent{\it#1\unskip.\kern.4em}}
\def\enddf{\medbreak}
\let\rk\df\let\endrk\enddf
\let\Roster\bgroup\let\endRoster\egroup
\def\\{}\def\text#1{\hbox{\rm #1}}
\def\mop#1{\mathop{\rm\vphantom{x}#1}\nolimits}
\def\sam{}\def\prose{}\def\MaxReferenceTag#1{}
\def\qedbox{\vrule width2mm height2mm\hglue1mm\relax}
\def\qed{\ifmmode\qedbox\else\hglue5mm\unskip\hfill\qedbox\medbreak\fi\rm}

\def\Smallfonts{}

\let\Item\item
\def\cite#1{{\bf[#1]}}
\def\Em#1{{\it #1\/}}\let\em\Em
\def\Bib#1\par{\bigbreak\bgroup\centerline{#1}\medbreak\parindent30pt
 \parskip2pt\frenchspacing\par}
\def\endBib{\par\egroup}
\newdimen\Overhang
\def\rf#1{\par\noindent\hangafter1\hangindent=\parindent
     \setbox0=\hbox{[#1]}\Overhang\wd0\advance\Overhang.4em\relax
     \ifdim\Overhang>\hangindent\else\Overhang\hangindent\fi
     \hbox to \Overhang{\box0\hss}\ignorespaces}

\def\bbN{{\Bbb N}}
\def\bbZ{{\Bbb Z}}
\def\Coordinates{\bigbreak\bgroup\parindent=0pt\obeylines}
\def\endCoordinates{\egroup}

\long\def\omit#1{} 
\newcount\NRcount\NRcount0
\def\NR{\advance\NRcount1\Nr{\the\NRcount.}}
\def\IT{\advance\NRcount1\Item{\the\NRcount.}}

\def\SetRev{\def\track{fellow travel}
\def\tracks{fellow travels}
\def\tracked{fellow travelled}
\def\tracking{fellow traveller}}

\SetRev
\def\bull{\item{\raise1pt\hbox{\sam$\bullet$}}}

\def\Im{\mop{Im}}

\def\int{\mop{int}}

\def\Autstruct{\mop{S\frak A}}
\def\AAutstruct{\mop{\frak A}} 
\def\BAutstruct{\mop{BS\frak A}}

\def\BAAutstruct{\mop{B\frak A}} 
\def\len{\mop{len}}

\def\cl{\mop{cl}}

\def\extY{\widehat{\cal Y}}
\def\vert{\mop{vert}}
\def\edge{\mop{edge}}
\let\eval\Overline
\def\init{\partial_0}\def\term{\partial_1}
\def\\{}
\def\Title{Automatic Structures and Boundaries\\ for Graphs of Groups}
\def\Author{Walter D. Neumann and Michael Shapiro}
\SetTitleHead{\Title}
\SetAuthorHead{\Author}

\H \Title
 
\Author\footnote{*}{Both authors
acknowledge support from the NSF for this research.}
 
{\Smallfonts\narrower
\HHH Abstract
 
We study the synchronous and asynchronous automatic structures on the
fundamental group of a graph of groups in which each edge group is
finite.  Up to a natural equivalence relation, the set of biautomatic
structures on such a graph product bijects to the product of the sets
of biautomatic structures on the vertex groups.  The set of automatic
structures is much richer.  Indeed, it is dense in the infinite
product of the sets of automatic structures of all conjugates of the
vertex groups. We classify these structures by a class of labelled
graphs which ``mimic" the underlying graph of the graph of groups.
Analogous statements hold for asynchronous automatic structures. We
also discuss the boundaries of these structures.\par}

\HH 1. Introduction

Given a group $G$, there is a natural equivalence relation on the set
of synchronous or asynchronous automatic structures on $G$.  Namely,
two such structures, $L$ and $L'$ are equivalent (written $L\sim L'$)
if there is a constant $K$ so that whenever a word of $L$ and a word
of $L'$ represent the same element of $G$, these two words
asynchronously $K$-fellow travel each other.  (For definitions, see
below.)  This leads \cite{NS1} to introduce  $\Autstruct(G)$,
$\BAutstruct(G)$, $\AAutstruct(G)$ and $\BAAutstruct(G)$, the sets of
(respectively) automatic, biautomatic, asynchronously automatic, and
asynchronously biautomatic structures on $G$ up to equivalence.

Currently, information is fairly scarce about these sets. 
$\AAutstruct(G)$ has been computed for $G$ virtually abelian,
virtually free (\cite{NS1}), or virtually a surface group (\cite{B}).
In the latter two cases it is a single point.
$\Autstruct(G)=\BAutstruct(G)$ and is a single point if $G$ is word
hyperbolic. $\BAutstruct(G)$ has been computed, and $\Autstruct(G)$ is
fairly well understood when $G$ is a geometrically finite hyperbolic
group (\cite{NS2}).  By contrast, $\AAutstruct(G)$ is very large and
poorly understood when $G$ is the fundamental group of a closed
hyperbolic 3-manifold group which fibers over the circle.  Notice that
here $G$ is an HNN-extension of a hyperbolic surface group.  However,
its unique automatic structure does not arise from the automatic
structure on the surface group.  Indeed, the surface group is not
rational in this automatic structure.

In this paper, we will study these sets when $G=\pi_1(\cal Y)$ is the
fundamental group of a finite graph of groups $\cal Y$ in which each
edge group is finite.  The assumption of finite edge groups turns out
to ensure that each conjugate of a vertex group is rational with
respect to any (synchronous or asynchronous) automatic structure on
$G$.  Consequently, there is a natural map
$$\AAutstruct(G)\to\prod_{H \in \cal H} \AAutstruct(H),$$
where $\cal H$ denotes the set of conjugates of vertex groups.  One of
our main results is that this map is injective with dense image
(Theorem 3.9).  Moreover, $\Autstruct(G)$ is just the inverse image of
$\prod_{H\in\cal H} \Autstruct(H)$ under this map.  In \cite{NS2}
there is an analogous injection for synchronous automatic structures
on a geometrically finite hyperbolic group with maximal parabolics
playing the role that vertex groups play here.

Now if $H \in \cal H$, we have $H=gG_Vg^{-1}$ where $G_V$ is a vertex
group.  Given such a $g$, it is natural to look for 
\omit{``$G_V$-free'' $h$} ``minimal'' $h$
so that $H=hG_Vh^{-1}$.  We shall see that if $H\ne G_V$ 
there is an $F_E$ orbit of
such values, where $F_E$ is the edge group associated to an edge $E$
incident at $V$.  Thus, to specify an asynchronous automatic structure
on $G$ it is only necessary to specify a choice of structure $[L_h]
\in \AAutstruct(G_V)$ for each such $h$, for this, in turn, specifies
a structure on $H=hG_Vh^{-1}$.  In a sense which we shall make clear,
the choice of $L_h$ must be equivariant with respect to the action of
$F_E$ on the set of such $h$.  This allows us to classify
$\AAutstruct(G)$ in terms of maps which we call \Em{regular
deployments} (Theorems 3.3 and 3.8).

There is a more concrete way to classify $\AAutstruct(G)$, and that is
in terms of objects which we call minimal special $\cal Y$-graphs.
Roughly, a $\cal Y$-graph $\cal X$ is a finite labelled graph which
maps onto the underlying graph of $\cal Y$.  Each vertex of $\cal X$
is labelled by an equivalence class of structures on the corresponding
vertex group of $\cal Y$.  Each edge from this vertex is labelled by a
rational subset of this vertex group.  The labelling must be
equivariant in terms of actions of the edge groups.  In the case of a
biautomatic or asynchronously biautomatic structure, a most efficient
$\cal Y$-graph is essentially the underlying graph of $Y$ with
biautomatic or asynchronously biautomatic structures at each vertex.
In particular, this gives bijections
$$\BAAutstruct (G) \to \prod_{V\in
\vert(\cal Y)}\BAAutstruct(G_V)$$
and
$$\BAutstruct (G) \to
\prod_{V\in \vert(\cal Y)} \BAutstruct(G_V),
$$ 
where $G_V$ denotes the vertex group at $V$.

As classifying objects these $\cal Y$-graphs have several
advantages over regular deployments.  They are finite objects and
they are easy to construct.  Unlike regular deployments, $\cal
Y$-graphs admit ``local'' modifications.  Finally, a $\cal Y$-graph
is easily turned into a generalized finite state automaton for the
structure which it determines.

In the final section of this paper we describe the boundary of an
asynchronous or synchronous automatic structure on $G=\pi_1(\cal Y)$. 
It is a ``tree completion'' of the disjoint union of the boundaries
for the automatic structures on the conjugates of the vertex groups.
The tree in question is the tree on which $G$ acts with quotient $\cal
Y$ (see for example \cite{Se}).

The assumption of finite edge groups in this paper may seem
restrictive.  However, as the hyperbolic 3-manifold example mentioned
above shows, it is necessary to ensure that the vertex groups are
rational.  In fact, even if one restricts to abelian edge groups, the
Heisenberg group is an example of the fundamental group of such a
graph of groups where the vertex group has plentiful automatic
structures but the group itself is not even asynchronously automatic.
Another example is $F_2\times\bbZ$, which can be seen either as
$F_2*_{F_2}$ or $\bbZ^2 *_\bbZ \bbZ^2$. In work in preparation we show
that $\Autstruct(F_2\times\bbZ)$ is quite large and unlikely to yield
to classification by techniques like the current ones.

\HH 2.  Background and definitions.

We start with a finitely generated group $G$ and a map from a finite
set $A=\{a_i\}$ into $G$ denoted by $a_i \mapsto \eval{a_i}$.  The set
of all finite strings $w=a_{i_1}\ldots a_{i_n}$ on elements of $A$
(including the empty string) forms a monoid under the operation of
concatenation.  We denote this monoid by $A^*$.  We define the
\Em{length} of $w=a_{i_1}\ldots a_{i_n}\in A^*$ to be $n$ and denote
this by $\len(w)$.  (The length of the empty word is $0$.)   The map
$a_i \mapsto \eval{a_i}$ extends to a unique monoid homomorphism from
$A^*$ to $G$ and we denote this extension by $w \mapsto \eval w$.  We
will assume that this map is onto.  We will also assume that $A$ is
supplied with an involution denoted by $a_i \mapsto a_i^{-1}$ and that
the evaluation map respects this, that is, $\eval{a_i^{-1}} = (\eval
{a_i})^{-1}$.  This allows us to form the \Em{Cayley graph} of $G$
with respect to $A$, $\Gamma=\Gamma_A$.  The vertices of $\Gamma_A$
are the elements of $G$.  There is a directed edge from $g$ to $g'$
labelled by $a\in A$ exactly when $g'=g\eval a$.  Thus there is
exactly an $A$'s worth of edges emanating from each vertex of $G$.
Since $A$ is finite, $\Gamma$ is locally finite.  Since $\eval A$
generates $G$, $\Gamma$ is path connected.  By making each edge of
$\Gamma$ isometric with the unit interval, we endow $G$ with a metric
$d_A=d$ called the \Em{word metric}.  That is, the distance between
two points of $\Gamma$ is defined to be the length of the shortest
path connecting them.  $G$ acts on $\Gamma$ by left translation, and
this action preserves distance. We take the \Em{length} of an element
of $G$ to be its distance from the identity, that is $\ell(g)=d(1,g)$.
A word $w\in A^*$ determines a path in $\Gamma$, which we also denote
by $w$, as follows. The path $w$ maps  the interval $[0,\len(w)]$
into $\Gamma$ by following at unit speed along the edge path in
$\Gamma$ based at $1$ and labelled by $w$.  We extend this to a
map of $[0, \infty)$ by setting $w(t)=\eval w$ for $t\ge \len(w)$. 

We call a subset of $A^*$ a \Em{language}.  A language $L$ is a
\Em{normal form} if $\eval L= G$.  Note that we do not require $L
\to G$ to be an injection.  We will say that a normal $L$ has
the \Em{asynchronous fellow traveller property} if there is a constant
$K$ so that given $w,w' \in L$ with $d(\eval w, \eval{w'})\le 1$,
there are monotone maps $\phi, \psi$ of $[0,\infty)$ onto itself so
that for all $t$, $d(w(\phi(t)),w'(\psi(t))) \le K$.  We say $L$ has
the \Em{synchronous fellow traveller property} if $\phi$ and $\psi$
can be chosen to be the identity.    

Given two normal forms $L, L'$, each with the asynchronous fellow
traveller property, we will say that they are \Em{equivalent} and
write $L\sim L'$ if $L\cup L'$ has the asynchronous fellow traveller
property.  We denote the equivalence class of $L$ by $[L]$.

Recall that a \Em{finite state automaton} $\cal A$ with alphabet $A$
is a finite directed graph on a vertex set $S$ (called the set of
\Em{states}) with each edge labelled by an element of $A$ and such
that different edges leaving a vertex always have different labels.
Moreover, a \Em{start state} $s_0\in S$ and a subset of \Em{accepted
states} $T \subset S$ are given.  A word $w\in A^*$ is in the
\Em{language $L$ accepted by $\cal A$} if and only if it defines a
path starting from $s_0$ and ending in an accept state in this graph.
We may assume there is no ``dead state'' in $S$ (a state not
accessible from $s_0$ or from which no accepted state is accessible).
Eliminating such states does not change the language $L$ accepted by
$\cal A$.

A language is \Em{regular} if it is accepted by some finite state
automaton.

We will also need the concept of a \Em{non-deterministic finite state
automaton}.  The difference is that a non-deterministic finite state
automaton is allowed to have several start states instead of just one,
different edges from a vertex may have the same label, and edges are
allowed to have empty label (such edges are called
\Em{$\epsilon$-transitions}).  A word is accepted by such an automaton
if it labels a path from a start state to an accept state.  This path
is allowed to traverse $\epsilon$-transitions.  It is a standard
result that the language of words accepted by a non-deterministic
finite state automaton is a regular language.

We shall also have occasion to use {generalized finite state
automata}.  A \Em{generalized finite state automaton} is defined just
like a non-deterministic finite state automaton except that the edges
are labelled by regular sublanguages of $A^*$ rather than by elements
of $A$.  This machine accepts a word $w$ if $w$ can be written as
$w_1\ldots w_k$ such that there is a corresponding directed edge path
$e_1\ldots e_k$ from a start state to an accept state such that $w_i$
is in the language labelling $e_i$ for each $i$.  It is a standard
fact that this language of accepted words is regular.

A \Em{(synchronous) automatic structure} for $G$ is a regular normal
form with the synchronous fellow traveller property.  It is a result
of \cite {ECHLPT} that every automatic structure has a sublanguage
which bijects to $G$.  Notice that if $L\subset L'$ and $L$ is an
automatic structure, then $L\sim L'$.  We will take the following  as
our definition of asynchronous automatic structure.  An
\Em{asynchronous automatic structure} for $G$ is a rational normal
form with the asynchronous fellow traveller property.  This is not
exactly equivalent to the use of the term in \cite{ECHLPT}.  Rather,
these are the \Em{non-deterministic asynchronous automatic structures}
of \cite{S1}.  Since every non-deterministic asynchronous automatic
structure contains an equivalent asynchronous automatic structure
which bijects to $G$, we will make no further distinction between the
two.  We will call an asynchronous automatic structure $L$
\Em{asynchronously biautomatic} if there is a constant $K$ so that if
$w, w' \in L$ with $\eval w = a\eval {w'}$ where $a \in \eval A \cup
\{1\}$, then there are reparameterizations $\phi$ and $\psi$ so that
for all $t$, $d(w(\phi(t)),aw'(\psi(t))) \le K$.  (Here $aw(\cdot)$ is
the translate of $w(\cdot)$ by $a$.)  We will call an automatic
structure
\Em{biautomatic} if there exists $K$ so that $\phi, \psi$ can be taken
to be the identity.  We take  $\Autstruct(G)$, $\BAutstruct(G)$
$\BAAutstruct(G)$ and $\BAAutstruct(G)$ to be respectively the sets of
automatic, biautomatic, asynchronously automatic, and
asynchronously biautomatic structures on $G$ up to equivalence.

Given an asynchronous or synchronous automatic structure $L$, we say
that $S\subset G$ is \Em{$L$-rational} if $\{w \in L : \eval w \in
S\}$ is regular.  It is a result of \cite{NS1} that $L$-rationality
depends only on the equivalence class of $L$.  Using the techniques of
\cite{GS}, one sees that if $H$ is an $L$-rational subgroup of $G$,
then $L$ induces an equivalence class of asynchronous respectively
synchronous automatic structures on $H$.

\rk Convention

We have pointed out that any asynchronous automatic structure contains
an equivalent one that bijects to $G$, so there is certainly no loss
of generality in assuming that all our structures are finite-to-one.
Since this simplifies some proofs, we will assume it from now on.
\endrk

\HH 3. Graphs of groups with finite edge groups

Let $\cal Y$ be a graph of groups.  We start by fixing notation.  The
underlying graph $Y$ of $\cal Y$ is a connected graph made of a finite
collection of vertices and a finite collection of unoriented edges.
We consider each unoriented edge as a pair of oriented edges and
denote the initial and terminal vertices of an oriented edge $E$ by
$\init E$ and $\term E$.  The reverse of an edge $E$ is denoted
$E^{-1}$.  To each vertex $V$ is associated a group $G_V$ and to each
edge $E$ is associated a group $F_E$ with $F_E=F_{E^{-1}}$.  Further,
to each edge $E$ is associated a pair of injections $\init:F_E\to
G_{\init E}$ and $\term:F_E\to G_{\term E}$, which are exchanged when
$E$ is replaced by $E^{-1}$.

Such a graph may be seen as instructions for building a group by
repeated free products with amalgamation and HNN-extensions. To do
this one takes a maximal tree $T\subset Y$.  Inductively one forms
free product with amalgamation for each edge of $T$.  One then
performs an HNN-extension for each edge not on $T$.  The resulting
group is determined up to isomorphism by the graph of groups. We refer
to it as the fundamental group of $\cal Y$, denoted $\pi_1(\cal Y)$.
For details see \cite{Se}, for example.

We describe a normal form for the elements of $G=\pi_1(\cal Y)$.  We
take the maximal tree $T$ to be fixed throughout.  We also choose a
fixed \Em{base vertex} $V_0\in \cal Y$.

\df Definition

For each edge $E$ of $Y$ we have an element $t_E\in G=\pi_1({\cal Y})$
as follows: $t_E$ is the stable letter associated to the edge $E$ if
$E$ is not in $T$ and $t_E=1$ if $E$ is in $T$.  In particular,
$t_{E^{-1}}=t_E^{-1}$. Then each element of $G$ can be written in the
\Em{normal form}
$$h=g_0t_{E_1}g_1\ldots t_{E_{m}}g_m\eqno{(*)}$$
where: 
\Roster
\Item{(1)} $E_1\ldots E_{m}$ is a path in $Y$ starting
at the base vertex $V_0$;
\Item{(2)} $g_0\in G_{V_0}$ and $g_i\in
G_{\term E_i}$ for $i=1,\ldots,m$;
\Item{(3)} if $E_{i+1}=E_i^{-1}$ then $g_i\notin\term(F_{E_i})$.  
\endRoster 
This expression is unique
up to the following two operations: 
\Roster
\Item{\bull} one can add or delete terminal words consisting of 
trivial $t_{E_i}$'s (subject to condition (3); such words are bounded
in length by the diameter of the maximal subtree $T$);
\Item{\bull} one can replace $g_it_{E_i}g_{i+1}$ by $(g_i\init(f))
t_{E_i}(\term(f^{-1})g_{i+1})$ for $f\in F_{E_i}$.
\endRoster\enddf

We shall assume from now on that all edge groups $F_E$ are finite.
Note that if $\cal Y$ includes an edge $E$ with $\term F_E=G_{\term
E}$ then this edge can be collapsed without changing $\pi_1(\cal Y)$
unless the edge is a loop, say $\init E=\term E = V$. In this case we
may also eliminate $E$ by replacing $G_V$ by $G_V \rtimes \bbZ$. This
$G_V \rtimes \bbZ$ is virtually cyclic so $\AAutstruct(G_V \rtimes
\bbZ)$ consists of a single point (cf.\ \cite{NS1}).  Thus, from the
point of view of computing asynchronous automatic structures on $G$ in
terms of asynchronous structures on the vertex groups, this
simplification of $\cal Y$ is harmless. If $\cal Y$ has no edge with
$\term F_E=G_{\term E}$ we will say $\cal Y$ is \Em{reduced}.

To simplify later notation we define an extended graph $\extY$ by
adding a new \Em{base edge} $E_0$ to $\cal Y$, going from a new vertex
(which we will never need to refer to) to the base vertex $V_0$.  We
put $F_{E_0}=\{1\}$. Using the above normal form we define for each
edge $E$ of $\extY$:
$$\cal G_E=
\{h=g_0t_{E_1}\ldots t_{E_m}g_m \text{ as in $(*)$}: m\ge0, E_m=E,
g_m\in\term(F_E)\}\prose.
$$ 
In particular, $\cal G_{E_0}=\{1\}$.  We
stress that the only role of the base edge $E_0$ is to support the
notation $\cal G_{E_0}$. We do not include the reverse edge
$E_0^{-1}$.

The significance of these sets is that, as we will discuss in detail
in section 5, the disjoint union $\coprod \cal G_E/F_E$ over
${E\in\edge\extY}$ is in one-one correspondence with the vertices
of the tree on which $G$ acts with quotient $\cal Y$ and vertex and
edge stabilizers given by the data of $\cal Y$.  If $\cal Y$ is
reduced then each conjugate of a vertex group stabilizes precisely one
vertex of this tree, so the disjoint union $\coprod \cal G_E/F_E$ is a
set of representatives for the conjugates of vertex groups. For our
present purposes we formulate this as the following lemma.

\th Lemma 3.1

Suppose $\cal Y$ is reduced. Let $H=gG_Vg^{-1}\in \cal H$.  Then $H$
determines an edge $E$ of $\extY$ with $\term E=V$ and $h\in \cal G_E$
so that $H=hG_Vh^{-1}$. If $E'$ with $\term E'=V$ and $h' \in \cal
G_{E'}$ also satisfy $h'G_V{h'}^{-1}=H$ then $E=E'$ and $h' \in
h\term F_E$.  \endth

\pf Proof

The edge $E$ and $h\in \cal G_E$ can be found from $g$ in the
following manner.  If $V=V_0$ and $g\in G_{V_0}$ we take $h=1\in\cal
G_{E_0}$. Otherwise, write $g$ in normal form $(*)$, and delete any
final portion of $g$ lying in $G_V$.  Call the resulting expression
$h$. If the last letter of $h$ lies in some $G_{V'}$, we take $E$ to
be the last edge in the path in $T$ from $V'$ to $V$.  If the last
letter of $h$ is the stable letter of an edge $E_1$ not in $T$, we
take $E=E_1$ if $\term E_1=V$ and otherwise we take $E$ to be the last
edge in the path in $T$ from $\term E_1$ to $V$.

Clearly $H=hG_Vh^{-1}$. The uniqueness statement about $E$ and $h$
follows by noting that $h$ (up to the action of $\term F_E$) is
visible as the first half of the normal form of any $hxh^{-1}\in
hG_Vh^{-1}$ with $x\notin \term F_E$.  \qed

\df Definition

Whether $\cal Y$ is reduced or not, we define a \Em{deployment} to be
a map
$$\psi:\coprod_{E\in \edge\extY}\cal G_E \to \coprod_{V\in
\vert\cal Y}\AAutstruct(G_V)
$$ 
with finite image taking $\cal G_E$ to $\AAutstruct(G_{\term E})$ for
each $E$ and with the following equivariance property: the restriction
of $\psi$ to $\cal G_E$ is $F_E$-equivariant in the sense that
$\psi(h)=f\psi(hf)$ for $f\in\term F_E$.\enddf

There is some ambiguity in the notation $\psi(h)$, for a group element
$h$ can be in more than one $\cal G_E$.  Thus we are implicitly
thinking of $h$ as an element of the disjoint union $\coprod \cal
G_E$.  The particular $\cal G_E$ intended should be clear from
context.

We wish to see that an arbitrary deployment determines an equivalence
class of (possibly non-regular) languages with the fellow traveller
property which map onto $G$.  We start by choosing a convenient
alphabet.

\df Definition

The above normal form for elements $G$ gives embeddings of the vertex
groups $G_V\subset G$.  We can take a generating set $A$ for $G$ which
is a union of generating sets for the vertex groups together with a
generator $t_E$ for each edge $E\notin T$. We denote by $A_V$ the
subset of all elements of $A$ evaluating into $G_V$.  We may choose
our generators such that for each edge $E$ of $\cal Y$ we have a
subset of $A$ which bijects to $\term F_E$.  By identifying any
duplicates, we can assume that each element of a group $\term
F_E\subset G$ is represented by exactly one letter of $A$.  We refer
to the subset of $A$ evaluating into $\term F_E$ as $A_E$.  In
particular, $A_E=A_{E^{-1}}$ for $E\in T$. Also, there is a unique
element $e\in A$ which evaluates to $1\in G$ and is in every $A_E$. We
call $A$ a \Em{convenient alphabet for $G$}.  \enddf

Let $\psi$ be a deployment.  For each $h\in \cal G_E$, $\psi(h)$ is a
class in $\AAutstruct(G_{\term E})$.  Choose a language $L_{\psi(h)}
\in\psi(h)$ for each $\psi(h)$.  We thus have $L_{\psi(h)}=
L_{\psi(h')}$ whenever $\psi(h)=\psi(h')$.  We assume our alphabet $A$
is convenient and each $G_V$-language is over the alphabet $A_V$.  We
take $t_E$ to be the empty word for $E\in T$.  Recall that $E_0$
denotes the ``base edge'' that we added to $\cal Y$ with $\term
E_0=V_0$.  We denote the unique element of $\cal G_{E_0}$ by $1$ and
define
$$\eqalign{
L_\psi = \{u_0 t_{E_1} \ldots t_{E_m}u_{m} :~& m\ge0,~E_1\ldots E_m
\text{ forms a path based at }V_0, \cr
&u_0 \in L_{\psi(1)}, u_i \in 
L_{\psi(\eval {u_0 t_{E_1}\ldots t_{E_i}})} \text { for } i \ge 1,\cr
&\text{if $E_{i+1}=E_i^{-1}$ then }\eval{u_i}\notin \term
F_E\}\prose.}
$$ 
In particular, $\eval {u_0} \in G_{V_0}, \eval {u_i}\in G_{\term
{E_i}}$ for $i\ge 1$.

\th Lemma 3.2

$L_\psi$ has the asynchronous fellow traveller property and is
determined up to equivalence by $\psi$.\endth

\pf Proof

We first show that $L_{\psi}$ has the asynchronous fellow traveller
property.  The fellow traveller constant will be $1+\max\{\delta_h\}$
where $\delta_h$ is a fellow traveller constant for $\hat L_{\psi(h)}
:=\bigcup_{f\in A_E} fL_{\psi(h\eval f)}$. This union is an automatic
structure by the equivariance of $\psi$.

So suppose $w,w' \in L_\psi$ with $\eval w = \eval{w'a}$, $a\in A$.
These words determine based edge paths $p$ and $p'$ in $\cal Y$ up to
terminal segments lying in $T$.

First suppose we can take $p=p'$ so that $w=u_0t_{E_1} \ldots
t_{E_m}u_m$ and $w'=\break u'_0t_{E_1} \ldots t_{E_m}u'_m$.  If $m=0$
then $w=u_0$ and $w'=u'_0$ both lie in $L_{\psi(1)}$ so they
fellow-travel.  Otherwise, $\eval{u_0}$ and $\eval{u'_0}$ differ at
most by an element of $\init F_{E_1}$.  Thus
$\eval{u_0}=\eval{u'_0f_0}$ with $f_0\in A_{E_1^{-1}}$. Hence, again,
$u_0$ and $u'_0$ asynchronously fellow travel with the fellow
traveller constant of $L_{\psi(1)}$.

More generally, $\eval {u_0t_{E_1}\ldots u_{i-1}t_{E_{i}}}=\eval
{u'_0t_{E_1}\ldots u'_{i-1}t_{E_{i}}g_i}$ with $g_i\in A_{E_{i}}$.  We
assume inductively that the two word have asynchronously
fellow-travelled to this point with fellow-traveller constant $\delta$
as above.  If we put $h=\eval{u'_0t_{E_1}\ldots u'_{i-1}t_{E_{i}}}$
then $u'_i\in L_{\psi(h)}\subset \hat L_{\psi(h)}$ and $g_iu_i\in
g_iL_{\psi(h\eval{g_i})}\subset\hat L_{\psi(h)}$. Also,
$\eval{g_iu_i}=\eval{u'_if_{i}}$ where $f_{i}\in A_{E_{i+1}^{-1}}$ is
such that $f_it_{E_{i+1}}=t_{E_{i+1}}g_{i+1}$.  Thus $g_iu_i$ and
$u'_i$ asynchronously fellow-travel with the fellow-traveller constant
of $\hat L_{\psi(h)}$. Hence $u_0t_{E_1}\ldots t_{E_i} u_i$ and
$u'_0t_{E_1}\ldots t_{E_i}u'_i$ asynchronously $\delta$-fellow-travel.
Thus, by induction, $w$ and $w'$ asynchronously
$\delta$-fellow-travel.

We must now examine the case where we cannot take $p=p'$.  In this
case, we can choose the paths $p$ and $p'$ so that (say) $p'$ is an
initial segment of $p$.  Then the previous case will apply to $w$ and
$w'a$.

Notice that the above argument also shows that $[L_\psi]$ did not
depend on the choices $\{L_{\psi(h)}\}$.  For if we are given choices
$\{L'_{\psi(h)}\}$ giving $L'_\psi$, we repeat the argument using
$\{L''_{\psi(h)}\}$ where $L''_{\psi(h)}=L_{\psi(h)}\cup L'_{\psi(h)}$
and observe that $L_\psi\cup L'_\psi \subset L''_\psi$.  \qed

\df Definition

We say that $\psi$ is a \Em{regular deployment} if $L_{\psi}$ is
regular, and hence an asynchronous automatic structure, for some
choice of languages in the classes in $\Im(\psi)$. We will see in the
proof of the following theorem that the word ``some'' in this
definition can be replaced by ``any''. \enddf

\th Theorem 3.3

The map $\psi\mapsto[L_\psi]$ gives a bijection $\{\text{regular
deployments}\}\to\AAutstruct(G)$.\endth

\pf Proof

We start by constructing the inverse map.  That is, we construct a
deployment $\psi_L$ from $[L]\in\AAutstruct(G)$.

\th Lemma 3.4 \cite{BGSS}

Let $L\subset A^*$ be a finite to one rational structure for a group
$G$.  Then for each $g \in G$ there are only finitely many $y \in A^*$
so that for some $x,z \in A^*$, $xyz \in L$ and $\eval y = g$. \endth

\pf Proof

We suppose not, and let $\cal A$ be a finite state automaton for the
language $L$.  We then have $x_1y_1z_1,x_2y_2z_2,\dots \in L$ with
$\eval{y_1}=\eval{y_2}=\dots$. Among these we can find infinitely many
$x'_1y'_1z'_1,x'_2y'_2z'_2,\dots$ so that each of $x'_1, x'_2, \dots$
labels a path from the start state of $\cal A$ to a common state of
$\cal A$. Among these we can find infinitely many $x''_1y''_1z''_1,
x''_2y''_2z''_2,\dots$ so that each of $x''_1y''_1, x''_2y''_2, \dots$
labels a path from the start state of $\cal A$ to a common state of
$\cal A$.  But then $x''_1y''_1z''_1,x''_1y''_2z''_1,\dots \in L$ with
$\eval{x''_1y''_1z''_1}=\eval{x''_1y''_2z''_1}=\dots$, contradicting
the assumption that $L$ is finite to one. \qed

\th Lemma 3.5

Let $A$ be a convenient alphabet for $G=\pi_1(\cal Y)$.  Given an
asynchronous automatic structure $L'$ on $G$, we can choose an
asynchronous automatic structure $L$ so that $L\sim L'$ and $L \subset
A^*$.  Moreover $L$ can be chosen so that if $y \in A^*$, $xyz \in L$
and $\eval y \in \term F_E$ for some edge $E$ of $Y$, then $y$ has the
form $e^mfe^n$ with $e,f\in A$ and $\eval e =1$. $L$ is synchronous if
$L'$ is.  At the possible cost of turning a synchronous structure into
an asynchronous one, we can assume $y$ is the single letter $f$.\endth

\pf Proof

It is an observation of \cite{NS1}, based on a result of \cite{ECHLPT},
that given arbitrary monoid generating sets $A$ and $B$ for a group
$G$ and an synchronous or asynchronous automatic structure $L' \subset
B^*$ for $G$, there is a synchronous respectively asynchronous
automatic structure $L'' \subset A^*$ with $L'' \sim L'$.  Choose $A$
as in the lemma and $L''$ as just described. By the previous lemma,
there are only finitely many $y \in A^*$ so that $xyz \in L''$ and
$\eval y$ is in some $\term F_E$.  We consider such words that are not
already in $A\{e\}^*$.  For each such word $y$ there is a unique word
$w_y\in A\{e\}^*$ with the same length and value.

Let $\cal A''$ be a finite state automaton for the language $L''$.
Considering it as a finite graph and replacing it by a finite cover if
necessary, we may assume that every path in $\cal A''$ labelled by one
of these words $y$ is embedded. We now construct a nondeterministic
machine as follows: wherever we see a path labelled by one of the
words $y$ we add a new path from the beginning point of this path to
its end point labelled by $w_y$. Call the language of this machine
$N$.  The language we seek is obtained from the language of $N$ by
removing the regular language of all words containing one of the words
$y$ as a subword and is hence regular.  Since two such subwords may be
adjacent, we can only ensure that a word with value $\eval f$ has the
form $e^mfe^n$.  If we wish a language in which the $e$'s do not
occur, we simply replace each edge of a machine labelled $e$ by an
$\epsilon$-transition. \qed

We choose $L$ as in Lemma 3.5.  Then each element of $L$ has a
decomposition
$$w = u_0t_{E_{1}} \ldots t_{E_{m}} u_m\eqno{(**)}$$ 
where

\bull $t_{E_i}$ is the empty word if $E_i \in T$; 

\bull $E_{1} \ldots E_{m}$ is a path in $Y$ starting at $V_0$;

\bull $u_i \in (A_{\init E_{i}})^*$ for $i<m$, $u_m \in (A_{\term
E_{m}})^*$.

\bull if $E_{i+1}=E_{i}^{-1}$ then $\eval{u_i}\notin F_{E_{i}}$.

This decomposition of $w$ is unique up to the following two
operations.  The placement of each $t_{E_i}$, $E_i \in T$, is
determined up to one of the finitely many subwords $e^mfe^n$ of Lemma
3.5.  Second, we can adjoin or delete a terminal reduced path of
$t_{E_i}$'s lying in $T$.  We call this decomposition an \Em{edge path
decomposition} of $w$.

For each $h \in \cal G_E$, we define
$$\eqalign{N_h = \{ v \in A^* : \exists w& \in L\text{ edge path
decomposed as in $(**)$ above such that}\cr & \text{for some $i$, }
E_i=E, v=u_i, \eval{u_0t_{E_1}\ldots u_{i-1}t_{E_i}}=h\}\prose,\cr}
$$ 
and set
$$L_h=\bigcup_{f\in A_E}fN_{h\eval f}\prose.$$ 
(As usual we think of $h$ as lying in $\coprod \cal G_E$.  This saves
wear and tear on subscripts. In the next section it will be helpful to
write $L_{E,h}$ instead.)

\th Lemma 3.6

For each $h$ in $\cal G_E$, $L_h$ is an asynchronous automatic
structure for $G_{\term E}$.  There are finitely many distinct
languages $L_h$.  The equivalence class of $L_h$ depends only on $h$
and the equivalence class of $L$. If $L$ is a synchronous automatic
structure, then so is $L_h$. The assignment $h\mapsto [L_h]$ is
$F_E$-equivariant in the following sense: for $f\in F_E$ we have
$[L_{h}]=f[L_{hf}]$. In particular this assignment is a deployment
$\psi_L$ which depends only on $[L]$. \endth

\pf Proof

Assume that $h\in\cal G_E$.  We first check that $N_h$ is regular.  We
will express $N_h$ as the union of two sublanguages, determined by
whether the $E_{i+1}$ in the definition of $N_h$ equals $E^{-1}$ or
not, and build a nondeterministic finite state automaton for each of
these languages (the nondeterminism will consist only in possibly
having several start states).

The edge path decomposition of $L$-words induces similar
decompositions for subwords of $L$-words, which we will use in the
following.

Let $\cal A$ be a finite state automaton for $L$.  There are only
finitely many words in the prefix-closure of $L$ which evaluate to
$h$.  We take those which have an edge path decomposition ending with
$t_E$. We let $S_h$ be the collection of states of $\cal A$ reached by
these words.  Let $R'_E$ be the collection of states of $\cal A$ which
are accept states or have a path to an accept state of $\cal A$
labelled by a word that is not in $(A_{\term E})^*$ and has edge path
decomposition (as just described for $L$-subwords) $t_{E'}\ldots$ with
$\init E'=\term E$ and $E'\ne E^{-1}$.  Let $\cal A'_h$ be the
nondeterministic machine obtained from $\cal A$ by making $S_h$ the
set of start states, $R'_E$ the set of accept states, and deleting all
arrows labelled by letters not in $A_{\term E}$.  Let $N'_h$ be the
language accepted by the machine $\cal A_h$.  Let $R''_E$ be the
collection of states of $\cal A$ which have a path to an accept state
of $\cal A$ labelled by a word that is not in $(A_{\term E})^*$ and
has edge path decomposition $t_{E^{-1}}\ldots$. Let $\cal A''_h$ be
defined like $\cal A'_h$ but using $R''_E$ instead of $R'_E$ and let
$N''_h$ be the corresponding language.  Then $N_h=N'_h\cup(N''_h -
\{e\}^*A_E\{e\}^*)$. It is thus a regular language.  Moreover, it is
determined by the $E$ and the subset $S_h$ of the states of $\cal A$,
so there are only finitely many different languages $N_h$.

The language $L_h$ is now a finite union of regular languages, hence
regular.  Moreover, it is determined by $E$ and the family of subsets
$S_{h\eval f}, f\in A_E$, of the states of $\cal A$, so there are a
finite number of these languages.

We show that $L_h$ surjects onto $G_{\term E}$. For $g\in G_{\term E}$
let $w\in L$ be a word with value $hg$.  If $v$ is the largest
terminal segment of $w$ which lies in $(A_{\term E})^*$, then $w$
decomposes as $ut_Ev$ with $\eval{ut_E}=h\eval f$ and $\eval{v}=\eval
{f^{-1}}g$ with $f\in A_E$.  Then $fv\in fN_{h\eval f}\subset L_h$ and
$\eval{fv}=g$.

We next show $L_h$ has the asynchronous fellow traveller property. If
not, we could find $f_iv_i\in f_iN_{h\eval{f_i}}\subset L_h$,
$f'_iv'_i\in f'_iN_{h\eval{f'_i}}\subset L_h$, for $i=1,2,\ldots$, so
that $d(\eval{f_iv_i},\eval{f'_iv'_i})$ is bounded, but for any $K$,
there is some $i$ so that $f_iv_i$ and $f'_iv'_i$ do not
asynchronously $K$-fellow travel.  We would then have $u_1v_1w_1$,
$u'_1v'_1w'_1$, $u_2v_2w_2$, $u'_2v'_2w'_2, \dots\in L$ with $h\eval
f_i =\eval {u_i}$, $h\eval{f'_i}=\eval{u'_i}$ for each $i$.  We can
replace each $w_i$ and $w'_i$ with $x_i$ and $x'_i$ of bounded length.
Then for each $i$, $u_iv_ix_i, u'_iv'_ix'_i \in L$,
$d(\eval{u_iv_ix_i},\eval{u'_i,v'_ix'_i}) =
d(h\eval{f_iv_ix_i},h\eval{f'_iv'_ix'_i}) =
d(\eval{f_iv_ix_i},\eval{f'_iv'_ix'_i})$ is bounded, and yet there is
no $K$ so that each of the pairs $u_iv_ix_i$ and $u'_iv'_ix'_i$
asynchronously $K$-fellow travel.  This contradicts the assumption
that $L$ is an asynchronous automatic structure.

The same argument shows $L_h$ is synchronous if $L$ is.  If $L\sim L'$
we can apply the argument to $L\cup L'$ to see that, for a fixed $h$,
$[L_h]$ depends only on $[L]$.  Finally, the equivariance property is
immediate from the definition of $L_h$.\qed

\pf Proof of Theorem 3.3 continued

It is clear that $\psi\mapsto[L_\psi]$ maps $\{\text{regular}$\break
$\text{deployments}\} \rightarrow \AAutstruct(G)$.  We shall show
below that $[L]\mapsto \psi_{[L]}$ maps $\AAutstruct(G)\to
\{\text{regular}$ $\text{deployments}\}$.  Given this, it is easy to
see that these maps are mutual inverses.

For suppose we start with $[L]\in\AAutstruct(G)$. For $h\in\cal G_E$
we take $\hat L_h =N_h\cup\bigcup_{L_h \sim L_{h'}} L_{h'}$.  Since
this is a finite union of equivalent asynchronous automatic structures
on $G_{\term E}$, it is itself an asynchronous automatic structure.
We use the languages $\hat L_h\in\psi_{[L]}(h)$ to define
$L_{\psi_{[L]}}$. Then $L_{\psi_{[L]}}$ contains the language $L$ so
certainly $[L]= [L_{\psi_{[L]}}]$.

Similarly, if $\psi$ is a regular deployment then
$\psi_{L_\psi}(h)=[(L_\psi)_h]$ and $(L_\psi)_h$ contains the language
$eL_{\psi(h)}$.  Thus $\psi_{[L_\psi]}(h)= [L_{\psi(h)}]$ for each
$h$.  That is, $\psi_{[L_\psi]}(h)=\psi(h)$ for all $h$, so
$\psi_{[L_\psi]}=\psi$.

So, to complete the proof of 3.3, we need only prove that if $L$ is an
asynchronous automatic structure on $G$ then the deployment
$\psi_{[L]}$ is regular.  Let $\psi=\psi_{[L]}$.  We shall describe a
nondeterministic finite state automaton $\cal T$ for$L_\psi$, thus
showing $L_\psi$ is a regular language, so $\psi$ is a regular
deployment.

We assume $L$ is as in Lemma 3.5. The following is our key lemma.

\th Lemma 3.7

Given $K>0$, there exists a finite state automaton $\cal S=\cal S_K$
with the following properties:
\Roster
\Item{1.} It accepts any word $w\in A^*$;
\Item{2.} Suppose $w$ is a
word with value $h$ which asynchronously $K$-fellow travels a word of
$L$ with the same value.  Then the final state reached by $w$ in the
machine $\cal S$ tells one for each edge $E$ of $\extY$ whether $h \in
\cal G_E$, and if so, what the corresponding language $L_h$ is.
\endRoster\noindent In particular, if $L'$ is the prefix closure of
$L$, then the language $\{w\in L' : L_{\eval w}=L_h\}$ is regular for
any $h\in \coprod \cal G_E$.  \endth

\pf Proof

Recall that the language $L_h$ is determined by $E$ and the map
$f\mapsto S_{h\eval f}$ of $A_E$ to the power set of the set of states
of $\cal A$.  Moreover, at any point along a path $w$, $S_h$ is the
set of $\cal A$ states reached by paths in $\cal A$ having the same
value $h$ as our path's current value and decomposing as $ut_E$.  Thus
it behooves us to modify our machine $\cal A$ to make ``visible", the
invisible $t_E$'s when $E\in T$.  So suppose we take the alphabet
$A\cup\{r_E: E\in T\}$, and let $s_E=t_E$ if $E\notin T$, $s_E=r_E$ if
$E\in T$.  Then the reader can check that the language
$$\{u_0s_{E_1}
\ldots s_{E_m}u_m : w=u_0t_{E_1} \ldots t_{E_m}u_m \text{ is an edge
path decomposition of } w\in L\}
$$
is regular.  (Here is a sketch proof.  Add a loop with label $r_E$ to
every vertex of $\cal A$ for every $E\in T$.  This machine accepts the
language obtained from $L$ by adding arbitrary subwords in these
$r_E$'s.  The desired language is obtained from this by deleting any
word that has one of a certain finite collection of prohibited
subwords; it is hence regular.)

We take a deterministic machine for this language and replace each
$r_E$ edge with a $t_E$ edge which we take as an $\epsilon$
transition.  This is a machine for $L$ and we assume $\cal A$ is of
this form.

Now suppose $w$ $K$-fellow travels some path of $L$ with the same
value. Then any word of $L$ with the same value $(K+K_L)$-fellow
travels $w$, where $K_L$ is the fellow traveller constant for $L$.  It
thus suffices to keep track at each step along $w$ of what $\cal A$
states have been reached by paths which $(K+K_L)$-fellow travelled
ours and have final value in a $K$-neighborhood of our current value,
and, when one of these has the same value as $w$, whether its final
edge in $\cal A$ is a $t_E$ edge.  That is, the information we must
keep track of is an element of $\mop{Maps}(B,\cal P(S\times E(\cal
A)))$, where $B$ is a ball of radius $K+K_L$ in the Cayley graph, $S$
and $E(\cal A)$ are the sets of states and edges of $\cal A$
respectively, and $\cal P(\cdot)$ denotes the power set.

More precisely, to keep track of the desired information we use a
finite state automaton $\cal S$ with $\mop{Maps}(B,\cal P(S\times
E(\cal A)))$ as set of states.  For $\alpha,\beta\in\mop{Maps}(B,S
\times\cal P(E(\cal A)))$ and $a$ in our alphabet $A$, $\cal S$ has an
edge labelled $a$ from $\alpha$ to $\beta$ if and only if each
$\beta(g)$ consists of the set of pairs $(s,e)$ such that there exists
a path in $B\cup aB$ from a point $g_1\in B$ to $ag$ labelling a path
in $\cal A$ from a state $s'$ with $(s',e')$ in $\alpha(g_1)$ to $s$
with final edge $e$.  As start state we take the element $\sigma$ with
$\sigma(g)$ equal to the set of pairs $(s,e)$ so that $s$ reachable in
$\cal A$ from the start state by paths with value $g$ and final edge
$e$.  Any word $w\in A^*$ then defines a path in $\cal S$ from the
start state.  If $w$ asynchronously $K$-fellow travels an element of
$L$ with the same value then this path ends in a state $\alpha$ with
$p(\alpha(\eval f))=S_{\eval {wf}}$ for each $f\in A_E$ where $p$
denotes projection onto the first factor.  This state $\alpha$ thus
gives the desired information.\qed

We now return to the proof that $L_\psi$ is regular for $\psi=\psi_L$.
We choose regular languages $L_{\psi(h)}$ for each $\psi(h)$ to define
the language $L_\psi$.  By proving that this $L_\psi$ is regular, we
will also have proved the remark preceding Theorem 3.3.

We have already shown that $L$ asynchronously fellow travels one
choice of $L_\psi$ (namely, the one with $L_{\psi(h)}=\hat L_h$).
Hence, by Lemma 3.2 it asynchronously fellow travels any choice of
$L_\psi$.  Let $K$ be the fellow traveller constant for our particular
choice.  Let $\cal S$ be the machine of the above lemma.

Let $\cal A'$ be the disjoint union of machines for the languages
$L_{\psi(h)}$.  We shall construct a nondeterministic machine for the
language $L_\psi$ by adding some arrows to the product machine of
$\cal A'$ and $\cal S$. Namely, for each state $s$ of this product
machine and each edge $E$ of $\cal Y$ we will add an arrow labelled
$t_E$ from $s$ to the following state $t$, if it exists. The $\cal S$
component of $t$ is the one determined by the $t_E$-transition from
the $\cal S$-component of $s$.  The $\cal A'$ component is the start
state of the machine for $L_{h\eval{t_E}}$, where
$\psi(h\eval{t_E})=[L_{h\eval{t_E}}]$ is determined by $E$ and the
$\cal S$ component of $s$ as in the above lemma.  It is easy to see
that this machine performs as advertised.\qed

There is an entirely analogous version of Theorem 3.3 for synchronous
automatic structures.

\th Theorem 3.8

The bijection of Theorem 3.3 restricts to a bijection between
$\Autstruct(G)$ and the set of regular deployments whose images lie in
$\coprod_{V\in \vert\cal Y} \Autstruct(G_V)$.\endth

\pf Proof

Given a synchronously automatic structure $L$ on $G$, Lemma 3.6
ensures that the deployment $\psi_L$ takes its image in $\coprod_{V\in
\vert\cal Y}\Autstruct(G_V)$. It remains to check that if $\psi$
is a regular deployment whose image lies in $\coprod_{V\in
\vert\cal Y}\Autstruct(G_V)$ then we can find a synchronously
automatic structure $L \sim L_\psi$.

So suppose $\psi$ is such a deployment.  We choose synchronous
automatic structures with uniqueness $L_{\psi(h)}$ for each class in
$\Im \psi$. For each such language and each edge $E'$, if the language
occurs as $L_{\psi(h)}$ with $h\in\cal G_E$ and $\term E=\init E'$, we
choose a set of $F_{E'}$-coset representatives in the language
$L_{\psi(h)}$.  To do this we order our alphabet $A$.  This induces a
total order $\prec$ on $A^*$ by ordering first on length and then by
lexicographic order for words of a given length. Let
$$L'_{\psi(h),E'}
= \{ w \in L_{\psi(h)} : w \text{ is $\prec$-minimal among $w$ with }
\eval w \in \eval w F_{E'} \}\prose.
$$
By \cite{BGSS}, each of these languages is regular.  We take
$$\eqalign {L= \{ u_0t_{E_1} \ldots
t_{E_m} u_m : &\eval{u_0}\in G_{V_0}, \eval{u_i}\in G_{\term
E_i}\text{ for } i>0, \cr & u_i \in L'_{\psi (\eval{u_0t_{E_1} \ldots
t_{E_i}}),E_{i+1}} \text{ for } i < m, \cr &u_m \in L_{\psi(\eval{u_1
\ldots t_{E_{m}}})}\} \prose .\cr}
$$ 

Notice that by the normal form for graphs of groups, $L$ bijects to
$G$.  Further, $L \subset L_\psi$, hence $L\sim L_\psi$.  To see that
$L$ is regular, one builds a product machine based on $\cal S$ and
$\cal A'$, where here $\cal A'$ is the disjoint union of machines for
the languages $L'_{\psi(h),E'}$, and modifies this to accept
$L_{\psi(h)}$ in the final factor.  Finally we wish to see that $L$
has the synchronous fellow traveller property.  We repeat the argument
that $L_\psi$ has the asynchronous fellow traveller property, but with
the following observation.  Suppose $w={u_0t_{E_1} \ldots t_{E_m}
u_m}$ and $w' = {u'_0t_{E'_1} \ldots t_{E'_{m'}}u'_{m'}}$ with $w,w'
\in L$ and $\eval w = \eval{w'a}$.  Let $p=E_1\ldots E_m$,
$p'=E'_1\ldots E'_{m'}$.  If $p=p'$, then by our choice of coset
representatives, $u_0t_{E_0}\ldots t_{E_m}=u'_0t_{E'_0}\ldots
t_{E'_{m'}}$, and since $u_m$ and $u_{m'}$ synchronously fellow
travel, so, too, do $w$ and $w'$.  On the other hand, if we cannot
choose the decompositions so that $p=p'$, then (say) $p'$ is an
initial segment of $p$ and we in fact have $w'=wa$.\qed

\th Theorem 3.9

Let $L$ be an asynchronous automatic structure on $G=\pi_1(\cal Y)$.
Then each subgroup $H\subset G$ conjugate to a vertex group is
$L$-rational, and hence has an induced automatic structure $L_H$,
determined up to equivalence by $[L]$.  Let $\cal H$ be the set of all
conjugates of vertex groups. Then the map $[L]\mapsto
([L_H])_{H\in\cal H}$ defines maps
$$\AAutstruct(G)\to\prod_{H\in\cal H}\AAutstruct(H)\prose,\quad
\Autstruct(G)\to\prod_{H\in\cal H}\Autstruct(H)\prose,$$
which are injective and have dense image in the product
topology.\endth

\pf Proof

We replace $\cal Y$ by a reduced graph of groups (this is defined just
before Lemma 3.1). This does not affect any vertex with non-trivial
$\AAutstruct(G_V)$, and hence does not change the validity of the
theorem.  Let $L\subset A^*$ be an asynchronous automatic structure on
$G$ chosen according to Lemma 3.5.  Let $H=hG_Vh^{-1}$ with $h\in\cal
G_E$ be as in Lemma 3.1 and denote $V=\term E$.  By considering the
edge path decomposition one sees that any $w\in L$ with $\eval w\in H
- h(\term F_E)h^{-1}$ has the form $w=xyz$ with $\eval x,(\eval z
)^{-1}\in h\term F_E$ and $y\in (A_V)^*$.  Since there are just
finitely many possibilities for $x$ and $z$, the set of $L$-words of
this form is a regular sub-language of $L$.  Thus $H-h(\term
F_E)h^{-1}$ is rational, whence $H$ is, since $h(F_E)h^{-1}$ is
finite.

We also see that $\{w \in L : \eval w \in H\} \sim h L_h h^{-1}$.
Thus $[L_H]=h \psi_L(h) h^{-1}$.  In particular, the map $H\mapsto
[L_H]$ determines $\psi_L$ and hence $[L]$, so the map
$\AAutstruct(G)\to\prod_{H\in\cal H}\AAutstruct(H)$ is injective.  The
statement that this map has dense image is the statement that, if we
specify a structure $L_H$ for finitely many $H$, there is a structure
$L$ on $G$ which realizes these $L_H$'s. We defer the proof of this to
the next section.

The statement of the theorem in the synchronous case follows from the
above together with Theorem 3.8.\qed

\rk Remark

We can turn $L_h$ itself into an automatic structure on $H$ if we use
the evaluation map $A_V\to H$ given by $a\mapsto h\eval ah^{-1}$.
With this interpretation, $[L_h]=[L_H]$.\endrk

\def\e{{\bf e}}\def\v{{\bf v}}

\HH 4. $\cal Y$-graphs

Regular deployments are unsatisfactory as classifying objects.  This
is for several reasons. First is the fact that it is hard to specify
an arbitrary deployment and there is no convenient way to tell \em{a
priori} whether a deployment is in fact regular. This fact is
reflected strongly in the second, namely, that regularity is non-local
in following sense.  If one changes the value of a regular deployment
on one element $h\in\coprod\cal G_E$ one is likely to obtain a
non-regular deployment.

In this section we introduce a classifying object which avoids these
deficiencies. 

\df Definition

Given a finite graph of groups $\cal Y$ with finite edge groups, a
\Em{$\cal Y$-graph} $\cal X$ is a finite directed labelled graph $X$
with the following additional structure: 

\bull A map $\pi\colon X\to Y$ of underlying graphs is given.  A
vertex $\v$ of $\cal X$ with $\pi(\v)=V$ is called a $V$-vertex and
an edge $\e$ of $\cal X$ with $\pi(\e)=E$ is called an $E$-edge. This
is called the $\cal Y$-type of $\v$ or $\e$.

\bull A vertex $\v_0$ of $\cal X$ is chosen as \Em{start vertex} and
every vertex of $\cal X$ can be reached by a directed path from this
start vertex. (We may assume that $\v_0$ is a $V_0$-vertex, where
$V_0$ is the base vertex for $\cal Y$ chosen in Section 3.)

\bull Each $V$-vertex $\v$ is labelled by an element $[L_{\v}]\in
\AAutstruct(G_{V})$.

\bull Each edge $\e$ out of $\v$ is labelled by an $[L_{\v}]$-rational
subset $S_{\e}$ of $G_V$.  For each edge $E$ of $\cal Y$ out of $V$,
the labels on the $E$-edges out of $\v$ are disjoint.  Their union is
$G_V$ if $\v= \v_0$ or if $\v$ has an incoming edge of $\cal Y$-type
other than $E^{-1}$.  If $\v\ne \v_0$ and all incoming edges at $\v$
are $E^{-1}$-edges, their union is $G_V-\init F_E$.

\bull For each $V$-vertex $\v$ and each edge $E$ out of $V$, there is
a $F_{E}$-action on $\cal X$ which fixes all vertices except those
reached by one $E$-edge from $\v$.  This action respects labels in
the following sense.  For a vertex $\v'$ reached by an $E$-edge from
$\v$ we have $[L_{f\v'}]=f[L_{\v'}]$.  For an edge $\e$ departing
$\v$ with label $S_{\e}\subset G_{V}$ the edge $f\e$ has label
$S_{f\e}=S_{\e}f^{-1}$.  For an edge $\e$ departing a vertex
reached by an edge from $\v$ we have $S_{f\e}=fS_{\e}$.
\enddf

Given a $\cal Y$-graph, $\cal X$, the following choices determine a
language $L_{\cal X}$ for $G=\pi_1(\cal Y)$.  Choose a convenient
generating set $A$ for $G$. For each vertex $\v$ of $\cal X$,
choose an asynchronous automatic structure $L_{\v}\subset
(A_{\pi(\v)})^*$ in the class associated to $\v$.  For an edge $\e$
departing $\v$ let $L_{\e}$ be the sublanguage of words of $L_{\v}$
that represent elements of the set $S_{\e}$.  Let $T$ be a maximal
spanning tree in $\cal Y$.  For each edge $E$ of $\cal Y$ let $t_E$ be
as in Section 3, that is, it is the corresponding stable letter if
$E\notin T$ and the empty word if $E\in T$. Then $L_{\cal X}$ is
$$\eqalign{ L_{\cal X}=\{u_0t_{\pi(\e_1)}\ldots t_{\pi(\e_m)}u_m :~
&\e_1\ldots \e_m\text{ is a path in $\cal X$ from the start vertex,}\cr
&u_k\in L_{\e_{k+1}} \text{ for $k=0,\ldots ,k-1$,}\cr 
&u_m\in L_{\term \e_m}\prose,\cr 
&\eval{u_{k+1}} \notin\init(F_{\pi(\e_{k+1})})
\text{ if $\pi(\e_{k+1})=\pi(\e_k)^{-1}$} \}\cr}
$$

\th Theorem 4.1

The above language $L_{\cal X}$ is an asynchronous automatic structure on
$\pi_1({\cal Y})$ and depends, up to equivalence, only on the $\cal
Y$-graph $\cal X$.

Every asynchronous automatic structure on
$\pi_1({\cal Y})$  is equivalent to one constructed as above. \endth

\pf Proof

The proof that $L_{\cal X}$ has the asynchronous fellow traveller
property and is determined up to equivalence by $\cal X$ is just like
the proof of the analogous statement for a language determined by a
deployment (Lemma 3.2), and is left to the reader.  (In fact, it is
not hard to see that $L_{\cal X}$ is contained in a language
determined by the following deployment $\psi$. Given $h\in \cal G_E$,
we find a path $p$ in $\cal X$ from the start vertex of $\cal X$ whose
final edge $\e$ has $\pi(\e)=E$, and the language determined by $p$
contains a word with value $h$.  We take $\psi(h)= [L_{\term \e}]$.)

We must check that $L_{\cal X}$ is regular.  To do this it is helpful
to modify $L_{\cal X}$ by redefining $t_E$ for each $E\in T$
temporarily to be a new letter which evaluates to $1\in G$, rather
than the empty word. We first turn $\cal X$ into a generalized finite
state automaton $\cal A_{\cal X}$.  We do this by subdividing each
edge $\e$ of $\cal X$ into two edges.  We label the first of these by
$L_{\e}$ and the second by $t_{\pi(\e)}$.  The start state of $\cal
A_{\cal X}$ is the start vertex of $\cal X$.  We take all vertices of
$\cal A_{\cal X}$ to be accept states.  The language of this machine
contains $L_{\cal X}$.  In fact, $L_{\cal X}$ is exactly the
sublanguage of words containing no substring of the form $t_E u
t_E^{-1}$ with $u \in L_{\v}$, $\pi(\v) =\term E$, $\eval u \in \term
F_E$.  Since there are finitely many such strings, $L_{\cal X}$ is
regular as required.  If we now replace each letter $t_E$ with $E\in
T$ by the empty word we get our original $L_{\cal X}$ back, and it is
still regular.

We must check that every asynchronous automatic structure arises as
above.  Suppose $L$ is an asynchronous automatic structure on $G$.  We
assume that our language $L$ and alphabet $A$ are as in Lemma 3.5. We
shall construct a $\cal Y$-graph $\cal X$ for $L$ of a rather special
type.  The start vertex will have no incoming edges and each vertex
other than the start vertex will have incoming edges all of one $\cal
Y$-type.

Let $\extY$ be $\cal Y$ extended by a base edge as in Section 3, and
for $E$ an edge of $\extY$, let $\cal G_E$ be as defined in Section 3.
We refer to the edge path decomposition of elements of $L$ described
before Lemma 3.6.  For $h\in \cal G_E$ we define
$$
\eqalign{N^{E,h}=\{v \in A^* : \exists w&\in L\text{ with edge path
decomposition $w = u_0t_{E_{1}} \ldots t_{E_{m}} u_m$}\cr&\text{such
that for some $i$, } E_i=E, \eval{u_0t_{E_1}\ldots
u_{i-1}t_{E_i}}=h\cr &v=u_it_{E_{i+1}}\ldots t_{E_m}u_m\}\prose,\cr}
$$ 
and 
$$L^{E,h}=\bigcup_{\eval f\in\term f_E}fN^{E,h\eval f}\prose.
$$

We claim that for each $h \in \cal G_E$ the language $N^{E,h}$ is
regular. Let $\cal A$ be a machine for $L$. As in the proof of Lemma
3.6, we let $S_h$ be the set of states of $\cal A$ reached by words in
the prefix closure of $L$ which evaluate to $h$ and have an edge path
decomposition ending in $t_E$.  Then $N^{E,h}$ is the language
accepted by the machine obtained from $\cal A$ by making $S_h$ the set
of start states.  Thus $N^{E,h}$ is regular.  It is also determined by
the finite set of states $S_h$, so there are finitely many different
languages $N^{E,h}$. Thus there are also finitely many languages $L^{E,h}$,
and they are regular.   The language $L^{E,h}$ is determined by the map
$f\mapsto S_{h\eval f}$ of $A_E$ to the power set of the set of states
of $\cal A$.

The subset of $G$ onto which the language $L^{E,h}$ evaluates depends
only on $E$: it is the set of elements of $G$ whose normal form
decomposition with base vertex $\term E$ cannot start with
$t_{E^{-1}}$.  Note that, except for a finite number of elements of
$G$, if an element is distance $1$ from this set then it is also in
this set. It therefore makes sense to talk about the asynchronous
fellow traveller property for $L^{E,h}$, even though this language
does not surject to $G$.  We claim that $L^{E,h}$ has this property.
For suppose $fv\in fN^{E,h\eval f}$ and $f'v'\in f'N^{E,h\eval{f'}}$
with $\eval f, \eval {f'} \in \term F_E$ and
$d(\eval{fv},\eval{f'v'})\le 1$.  Then there exist $u,u'$ with $\eval
u=h\eval f$ and $\eval{u'}=h\eval{f'}$ so that $uv\in L$ and $u'v'\in
L$.  Since $uv, u'v'\in L$ and $d(\eval{uv},\eval{u'v'})=
d(\eval{fv},\eval{f'v'})\le1$, $uv$ and $u'v'$ asynchronously
fellow-travel. It follows that after reparameterization, $fv$ and
$f'v'$ also asynchronously fellow-travel.

For fixed $E$ and $h,h'\in\cal G_E$ it therefore also makes sense to
ask if $L^{E,h}\sim L^{E,h'}$.  We define
$$\widehat L^{E,h}=\bigcup_{L^{E,h'}\sim L^{E,h}}L^{E,h'}\prose.$$
Since this is a finite union of regular equivalent languages, it is
also regular with the asynchronous fellow-traveller property.

Each $\widehat L^{E,h}$ induces an asynchronous automatic structure
$\widehat L_{E,h}$ on $G_{\term E}$.\break  Namely, we define
$$\widehat L_{E,h}= \{u\in A^* :\exists v \in A^*, uv\in \widehat L^{E,h},
\eval u \in G_{\term E}\}\prose.$$ 
It is easy to see that this language is an asynchronous automatic
structure on $G_{\term E}$.  In fact, it is just
$\bigcup_{L^{E,h'}\sim L^{E,h}}L_{E,h'}$, where $L_{E,h'}$ is the
$L_{h'}$ of Lemma 3.6 (we are now making the edge $E$ explicit in our
notation).

We are now prepared to describe the $\cal Y$-graph, $\cal X$ of the
Theorem.  For each $E\in \edge\extY$ it has a vertex for each
$\widehat L^{E,h}$.  The vertex corresponding to $\widehat L^{E,h}$
projects to $\term E$ under $\pi$, and is labelled by $[\widehat
L_{E,h}]$. We take the vertex corresponding to $\widehat L^{E_0,1}$ to
be the start vertex.  Suppose that $E$ and $E'$ are edges of $\cal Y$ with
$\term E = \init E'$.  There is an edge $E'$ from the vertex for
$\widehat L^{E,h}$ to the vertex for $\widehat L^{E',h'}$ if the set
$$ S_{E'}= \{ g \in G_{\term E} : \widehat L^{E',h'} = 
\widehat L^{E',hg}\}$$
is not empty.  In this case $E'$ is labelled by $S_{E'}$ and
projects to $E'$ under $\pi$.

We must check that this defines a $\cal Y$-graph.  That is, we must
see that the set $S_{E'}$ is well defined, that it is $\widehat
L_{E,h}$-rational, that the labels on the $E'$ edges out of the vertex
for $\widehat L^{E,h}$ partition $G_{\term E}$ if $E'\ne E^{-1}$ and
partition $G_{\term E}-\term F_E$ if $E'=E^{-1}$, and that $\cal X$
has the appropriate equivariance properties.

So suppose $g \in G_{\term E}$.  We check that $\widehat L^{E',hg}$
depends only on $\widehat L^{E,h}$ and $g$, and does not depend on
$h$.  We first show that $L^{E',hg}$ depends only on $L^{E,h}$ and
$g$.  Now $w\in L^{E',hg}=\bigcup_{f\in A_{E'}}fN^{E',hg\eval f}$ if
and only if $w=fv$ with $f\in A_{E'}$ and there is $u$ with $\eval u =
hg \eval f$ and $uv \in L$, and $u$ ending with $t_{E'}$ in some edge
path decomposition of $uv$.  If there is such a $u$, it has the form
$u=xy$ with $\eval x = h \eval {f'}$, $\eval y = \eval{f'^{-1}}g
\eval f$, where $\eval {f'}$ in $\term F_E$.  Thus $w=fv \in L^{E', h
g}$ if and only if we find $f' y v \in L^{E,h}$ with
$\eval{f'y}=g\eval f$.  (The reader might want to draw a picture.)
Thus $L^{E',hg}$ depends only on $L^{E,h}$ and $g$, and not on $h$.  A
similar argument shows that $L^{E',hg} \sim L^{E',h'g}$ if $L^{E,h}
\sim L^{E,h'}$.  So by taking the appropriate unions we see that
$\widehat L^{E',hg}$ depends only on $\widehat L^{E,h}$ and $g$, and
thus $S_{E'}$ is well defined.

The fact that $S_{E'}$ is $\widehat L_{E,h}$-rational will follow
from the existence of the machine $\cal S$ of Lemma 3.7.  Recall that
that machine does the following: if $w$ is a word with value $h\in
\cal G_E$ which fellow travels some $L$-word with the same value then
the state of $\cal S$ reached by $w$ determines the language $L_{E,h}$
(called $L_h$ in Lemma 3.7).  The way it does this is by determining
the map $f\mapsto S_{h\eval f}$ of $A_E$ to the power set of the set
of states of $\cal A$.  But, as we saw above, this map also determines
$L^{E,h}$, and hence $\widehat L^{E,h}$.  Thus $L_{E,h}$ can be
replaced by $\widehat L^{E,h}$ in Lemma 3.7.

Now let $u$ be an $L$-word for $h$ and $v$ be a $\widehat
L_{E,h}$-word for $g$.  Then, by construction of $\widehat L_{E,h}$,
the word $w=uv$ fellow travels an $L$-word for $hg$.  We can thus test
$w$ with the machine $\cal S$ to see if $\widehat L^{E',hg}=\widehat
L^{E',h'}$.  That is, $g\in S_{E'}$ if and only if the word $v$ for
$g$ labels a path in $\cal S$ from the state reached by $u$ to a state
that determines the language $\widehat L^{E',h'}$.  This is a regular
condition, so $S_{E'}$ is $\widehat L^{E,h}$-rational, as required. 

The labels on the $E'$-edges out of the vertex for $\widehat L^{E,h}$
are disjoint, since $\widehat L^{E,h}$ and $g$ determine $\widehat
L^{E',hg}$, and their union is clearly $G_{\term E}-\term F_E$ if
$E'=E^{-1}$ and $G_{\term E}$ otherwise.

We must construct the $F_{E'}$ action and show that it respects
labels.  To this end, let $E''$ be an edge of $Y$ with $\term E'=\init
E''$.  We suppose that there is an edge $E''_\ell$ from the vertex for
$\widehat L^{E',h'}$ to the vertex for $\widehat L^{E'',h''}$.  If
$\eval f \in F_{E'}$, we let $f$ carry the vertex for $\widehat
L^{E',h'}=\widehat L^{E',hg}$ to the vertex $\widehat L^{E',hg\eval
{f^{-1}}}$.  As we have seen, this is well defined.  The label at this
vertex is $[L_{E',h'\eval{f^{-1}}}]$, which is $[fL_{E',h'}]$ as
required.  It now follows that the action of $\eval f$ on vertices
induces an action on edges out of the vertex for $\widehat L^{E,h}$,
and this action respects edge labels, for if $g\in S_{E'}$ then
$g\eval{f^{-1}}$ labels an edge from the vertex for $\widehat L^{E,h}$
to the $\eval f$ image of the vertex for $\widehat L^{E',h'}$.  In
particular, we have $S_{\eval f E'}=S_{E'}\eval {f^{-1}}$ as
required. In the same way, the $\eval f$ action on vertices induces
an action on edges whose initial vertex is moved by $\eval f$ carrying
(say) $E''_\ell$ to (say) $\eval fE''_\ell$ so that $S_{\eval
fE''_\ell} = \eval f S_{E''_\ell}$.  This completes the proof that
$\cal X$ is indeed a $\cal Y$-graph.

Finally, we must check that $L$ is equivalent to $L_{\cal X}$.  It is
an easy induction on free product length that $L$ and $L_{\cal X}$
determine the same deployment.  Alternatively, one may note that
$L_{\cal X}$ is equivalent to $L$ since it contains $L$ as a
sublanguage. \qed

\rk Remark

Call a $\cal Y$-graph $\cal X$ \Em{special} if it satisfies:

\bull the start vertex $\v_0$ has no incoming edges;

\bull each vertex $\v\ne \v_0$ has incoming edges of just one $\cal
Y$-type.

Then the $\cal Y$-graph $\cal X$ for $L$ constructed in the above
proof is special, and it is not hard to verify that it is minimal with
this property, in the sense that any other special $\cal Y$-graph
$\cal X'$ defining a language equivalent to $L$ can be mapped to $\cal
X$ by a graph mapping that respects the $\cal Y$-type of vertices and
edges, respects vertex labels, and also respects edge labels in the
sense that the rational set associated to an edge of $\cal X'$ is
contained in the rational set associated to corresponding edge of
$\cal X$.  Thus minimal special $\cal Y$-graphs actually classify
asynchronous automatic structures on $G=\pi_1(\cal Y)$.  They are,
however, not always efficient classifying objects, in that one can
often find a much smaller non-special $\cal Y$-graph to describe the
same structure, as we will now describe.\endrk

Let $A$ be a convenient alphabet for $G=\pi_1(\cal Y)$. Suppose $\cal
X$ is a $\cal Y$-graph and let $L_{\cal X}\subset A^*$ be, as in
Theorem 4.1, the language of words labelling paths in $\cal X$ from
the start vertex. For any vertex $\v$ of $\cal X$ we can define
similarly the language $L_{\v}$ of words in $A^*$ that label paths
starting at $\v$.  (If $\cal X$ is as constructed in the above proof
and $\v$ is the vertex corresponding to $\widehat L^{E,h}$ then
$L_\v\sim \widehat L^{E,h}$.)  Now suppose that for some vertex $V$ of
$\cal Y$ we have $V$-vertices $\v$ and $\v'$ of $\cal X$ such that
$L_{\v}\cup L_{\v'}$ has the asynchronous fellow-traveller property.
We can then attempt to create a smaller $\cal Y$-graph by identifying
the vertices $\v$ and $\v'$ of $\cal X$.  If $E$ is an incoming edge
at $V$ then we have a $F_E$-action on $\cal X$ which permutes the
$V$-vertices, so we must do this identification equivariantly.  There
is no guarantee that we can do this, for if $L_{\v}\cup L_{\v'}$ and
$L_{\v'}\cup L_{\v''}$ have the asynchronous fellow-traveller property,
we cannot deduce that $L_{\v}\cup L_{\v''}$ does.  However, if we can
do this identification equivariantly, we obtain a smaller $\cal
Y$-graph for $[L]$.  

Even if one can collapse $\cal X$ as above, there may be several
inequivalent ways of doing so.  Indeed, it is not hard to find an
example of a graph of groups $\cal Y$ with finite edge groups for
which the $\cal Y$-graph $\cal X$ constructed in the proof of Theorem
4.1 can be collapsed to several inequivalent ``minimal'' $\cal
Y$-graphs. 

There are some situations in which the above collapse is clearly
possible. For example, if one has a $\cal Y$-graph $\cal X$ with no
``special'' vertices --- that is, every vertex $\v$ has more than one
$\cal Y$-type of incoming edge --- then every $L_{\v}$ surjects to $G$
so the condition that $L_{\v}\cup L_{\v'}$ have the asynchronous

fellow traveller property defines an equivalence relation on the
vertices of $\cal X$. It is then easy to see that there is a $\cal
Y$-graph that can be obtained by collapsing $\cal X$ as above.

Another case is when the language $L$ is asynchronously biautomatic.
Each $L_{\v}$ is equivalent to a sublanguage of a translate of $L$.
But, by definition of biautomaticity, any translate of $L$ is
equivalent to $L$. We can thus collapse all $V$-vertices to a single
vertex for each $V$.  This collapse is clearly equivariant and extends
trivially to edges, so we see that there is a $\cal Y$-graph $\cal X$
whose underlying graph $X$ is isomorphic to the underlying graph $Y$
of $\cal Y$. This proves part of:

\th Theorem 4.2

$L$ is an asynchronously biautomatic structure on $\pi_1(\cal Y)$ if
and only if $L$ has a $\cal Y$-graph $\cal X$ for which $\pi$ is an
isomorphism of the underlying graphs of $\cal X$ and $\cal Y$ and the
structure at each vertex $v$ of $\cal X$ is an asynchronously
biautomatic structure for $G_{\pi(v)}$.

The corresponding statement holds also with ``asynchronously
biautomatic'' replaced by ``biautomatic''.\endth

\pf Proof

We first point out that a asynchronous or synchronous biautomatic
structure on $G$ induces asynchronous or synchronous biautomatic
structures on the vertex groups, since they are rational subgroups. 
Thus the labels on the vertices of the above $\cal Y$-graph $\cal X$
are as claimed.

Now suppose we have a $\cal Y$-graph $\cal X$ as in the theorem.  We
check that the ensuing structure $L$ is asynchronously biautomatic.
So suppose $A$ is a convenient alphabet we have chosen languages at
each vertex, and suppose also that $w=u_0t_{E_1} \ldots t_{E_m}u_m$ is
an edge path decomposition of a word in the resulting language.  We
must show that if $a \in A$, the word for $\eval{aw}$ asynchronously
fellow travels $w$.  There are several cases.  We use the same symbol
for a vertex of $\cal Y$ and the corresponding vertex of $\cal X$.

We suppose first that ${u_0}\in (A_{V_0})^*$ is not empty and $\eval a
\in G_{V_0}$.  We take $u'_0$ to be the word in $L_{V_0}$ for $\eval
{au_0}$.  Then $w'=u'_0t_{E_1} \ldots t_{E_m}u_m$ is an accepted word
for $\eval {aw}$.  Since $L_{V_0}$ is asynchronously biautomatic,
$w'$ and $\eval a w$ asynchronously fellow travel as required.  If
$L_{V_0}$ is biautomatic, they synchronously fellow travel.

We now suppose that $u_0$ is non-empty and that $\eval a \notin
G_{V_0)}$.  Then we have $w'=au_0t_{E_1} \ldots t_{E_m}u_m$ an
accepted word and again $\eval a w$ and $w'$ appropriately fellow
travel, unless it happens that $a=t_{E}$, $\eval u_0\in \term F_E$,
and the first letter after $u_0$ in $w$ is $t_{E_i}=t_{E^{-1}}$.  In
this case $w'= f' u_{i+1}\dots u_m$, for suitable $f'\in A_{E^-1}$, is
the word we seek, and again $\eval a w $ and $w'$ synchronously or
asynchronously fellow travel as required.

The remaining cases are similar and are left to the reader.\qed

\pf Proof of Theorem 3.9 (completed)

We need to show the map of Theorem 3.9 has dense image.
Suppose that every vertex group $G_V$ has an asynchronous automatic
structure. We will first show that $G=\pi_1(\cal Y)$ has at least one
asynchronous automatic structure (this follows from the methods of 
\cite{S2}, but we give a proof here for completeness). We shall need
the following lemma.

\th Lemma 4.3

Suppose that $E$ is an edge of $\cal Y$, with $\init E = V$, and
suppose we are given $[L_V] \in \AAutstruct(G_V)$. Then there is a
partition of $G_V$ into distinct $[L_V]$-rational sets $S_{\eval f}$,
$\eval f \in F_E$ so that $F_E$ acts on the right to permute these
sets $\{S_{\eval f}\}$.  \endth

\pf Proof

We assume $L_V \in [L_V]$ is a structure with uniqueness.  We take
$L'$ to be those words $w \in L_V$ that are least in dictionary order
among the words that evaluate into $\eval w F_E$.  Since $F_E$ is
finite, and it is easy to check dictionary order by means of a finite
state automaton, the language
$$L''=\{(u,v)\in L_V\times L_V : \eval u \in \eval v F_E \text{ 
and $u$ precedes $v$}\}
$$ 
is the language of an asynchronous two tape automaton.  It follows
that $L'=L_V - p_2(L'')$ is regular. (Here $p_2$ denotes projection
onto the second factor.) We take $S_1 = \eval L'$, and for each $\eval
f \in F_E$, we take $S_{\eval f}=S_1\eval f$.  It is easy to check
that each of these is $L_V$-rational. \qed

We now construct a special $\cal Y$-graph $\cal X'$.  This graph will
have a vertex for each edge $E$ of $\extY$ and each element of $\term
F_E$.  For fixed $E$ these vertices will constitute a $\term
F_E$-orbit, and they will be labelled by the orbit of structures
$[\term fL_V]\in \AAutstruct(G_V)$, $f\in F_E$.  The above lemma
allows us to put in edges in an equivariant fashion to complete the
$\cal Y$-graph $\cal X'$.  Let $[L']$ be the structure determined by
$\cal X'$.

Recall that $\cal H$ is the set of conjugates of vertex groups in
$G=\pi_1(\cal Y)$ and we are trying to  show that the map
$$L\mapsto(L_H)_{H\in\cal H}~~:\Autstruct(G)\to\prod_{H\in\cal
H}\Autstruct(H)$$ 
has dense image.  We must show that if $H_1,\ldots,H_n$ are distinct
groups in $\cal H$ and we are given $[L_{H_1}],\ldots,[L_{H_n}]$ in
$\AAutstruct (H_1), \ldots,\AAutstruct (H_n)$, we can find $[L]\in
\AAutstruct (G)$ so that $[L]$ induces $[L_{H_1}],\ldots,[L_{H_n}]$ on
$H_1,\ldots,H_n$.
We shall modify the structure $L'$ described above to do what is
required. 

As discussed at the beginning of the proof of Theorem 3.9, we may
assume that $\cal Y$ is reduced. For each $i=1,\ldots n$ choose $E_i$
and $h_i\in \cal G_{E_i}$ as in Lemma 3.1 with $H_i=h_iG_{\term
E_i}h_i^{-1}$.  For each $i$ the normal form for $h_i$ determines a
path in $\cal Y$ starting at the base vertex $V_0$.  Inclusion of
these paths in each other induces a partial order on the $h_i$ and
hence on the $H_i$. We may assume the ordering $H_1,\ldots,H_n$
respects this partial order.  We will describe a modification of $\cal
X'$ to make the structure on $H_i$ equal to the desired one without
changing the structure on any $H_j$ which is earlier in the partial
order. Repeating this iteratively for $i=1,\ldots,n$ then proves the
theorem.

Thus suppose $i$ is chosen and write $h=h_i$, $H=H_i$. By taking a
cover of $\cal Y$ if necessary, we may assume that the path $\sigma$
in $\cal Y$ determined by $h$ is embedded. There is an induced
covering of $\cal X'$ and we replace $\cal X'$ by this covering. We
choose a lift $\sigma'$ of the path $\sigma$ to $\cal X'$.

\def\w{{\bf w}} Suppose $\sigma'$ has length at least 2 (we leave the
case that it is shorter to the reader). Let the final two edges of
$\sigma'$ be $\e'$ and $\e$, so $\term \e'=\init \e =\w$ say. Let
$\pi(\e')=E'$, $\pi(\e)=E$, $\pi(\w)=W$. There is a word $u_0\ldots
t_{E'}ut_E$ labelling the path $\sigma'$ and evaluating to $h$. Then
$\eval u\in S_{\e}$.  We delete $\eval u$ from $S_{\e}$ and establish
a new edge out of $\v$ to a new vertex.  We label the new edge with
$\{\eval u\}$ and the new vertex with $[h^{-1}L_Hh]$.  We let the
edges out of this new vertex duplicate the edges out of $\term \e$. 
Likewise, for each $f\in \term F_E$ we delete $\eval uf$ from
$S_{f(\e)}$ and establish a new edge with label $\{\eval uf\}$ to a
new vertex labelled $[(hf)^{-1}L_Hhf]$.  For each of the vertices
$f(\w)$ with $f\in \term F_{E'}$ we perform the same operation,
constructing new edges to the vertices we have just added.  This
produces a new $\cal Y$-graph for a structure which induces the
desired structure on $H=hG_Vh^{-1}$ and has not changed the induced
structure on any earlier $H_j$.\qed

\rk  Remark

In the proof of Theorem 4.1, we were required to show that a $\cal
Y$-graph $\cal X$ determines a regular language $L_{\cal X}$, and to
do this, we turned $\cal X$ into a generalized finite state automaton
$\cal A_{\cal X}$ which almost accepted the language in question.  To
obtain the desired language, we only needed to delete those words
containing subwords of the form $t_E f t_E^{-1}$ where $\eval f \in
\term F_E$.  In fact, there is a straight forward procedure for
turning an $\cal Y$-graph $\cal X$ into a generalized finite state
automaton which accepts $L_{\cal X}$ itself.  The method here is to
build an generalized finite state automaton $\cal B_{\cal X}$ whose
underlying graph projects to that of $\cal A_{\cal X}$. For each
vertex of $v$ of $\cal X$, and each edge $E$ into $v$, there are two
vertices in $\cal B_{\cal X}$.  One of these is reached only by
elements of $\term F_E$, and there are no $\pi(E^{-1})$ edges out of
this vertex.  The other is reached by all elements not in $\term F_E$.
This latter has a full armamentarium of edges out of it.  The
interested reader may wish to fill in the details along the lines of
the proofs of Lemmas 1.1 and 3.1 in \cite{S2}.  \endrk

\def\cl{\mop{cl}}\def\TC#1{#1\hat{\,}}
\HH 5. The boundary

We recall the boundary of an asynchronous automatic structure, as
defined in \cite{NS1}. Let $L\subset A^*$ be an asynchronous automatic
structure on a group $G$. As usual, we assume $L$ is finite to one.
An \Em{$L$-ray} is an infinite word $w\in A^\bbN$, all of whose
initial segments are initial segments of $L$-words.  Two rays are
\Em{equivalent} if they asynchronously fellow travel (at a distance
that may depend on the rays). The \Em{boundary} of $L$ is the set
$\partial L$ of equivalence classes of rays with the following
topology.  For an $L$-rational subset $R$ of $G$, define $\partial R$
to be the set of rays which fellow travel $R$ (that is, travel in a
bounded neighborhood of $R$; the bound may depend on the ray).  These
sets form a basis of closed sets for a topology on $\partial L$. This
boundary can be attached to $G$: the sets $\cl R:=R\cup\partial R$ are
a basis of closed sets for a topology on $\cl G:=G\cup\partial L$
which has $\partial L$ as a closed subspace and $G$ as an open
discrete subspace.  (This topology is the ``rational topology'' of
\cite{NS1}. Other topologies on $\partial L$ are also discussed there.)

In \cite{NS1} the ``rehabilitated boundary'' $\widehat {\partial L}$ is
also discussed, which appears to be an appropriate notion for groups
with large abelian subgroups.  A subset $\sigma$ of $\partial L$ is
called an \Em{abstract simplex} if, for any choice of a neighborhood
in $\cl G$ for each point of $\sigma$, the intersection of these
neighborhoods is non-empty.  This makes $\partial L$ into a
topological abstract simplicial complex, the geometric realization of
which is the \Em{rehabilitated boundary} $\widehat{\partial L}$.

Let $B$ be a tree and $\TC B=B\cup C$ be its end
compactification.  Given a continuous map $p$ of a space $X$ to $B$, we
define the \Em{tree completion} of $X$ with respect to $p$ as the
disjoint union
$$\TC X=X\cup C\prose,$$
with the smallest topology for which $X$ is a subspace and the induced
map $\TC p\colon \TC X\to \TC B$ is continuous.  It is
an easy exercise to see that $\TC X$ is compact if and only if
$p$ is a proper map and is Hausdorf if and only if $X$ is Hausdorf.

Now let $\cal Y$ be a graph of groups with finite edge stabilizers and
$G=\pi_1(\cal Y)$.  By \cite{Se}, there is a $G$-tree $B$ with $B/G$
equal to the underlying graph $Y$ of $\cal Y$ and with edge and vertex
stabilizers given by the data of $\cal Y$.  Let $\pi\colon B\to Y$ be
the projection. The stabilizer of a vertex $v$ of $B$ is therefore a
conjugate $H_v$ of the vertex group $G_{\pi(v)}$.

As in section 3, for $[L]\in \AAutstruct(G)$ and $H$ a conjugate of
a vertex group, $[L_H]$ denotes the induced structure on $H$.

\th Theorem 5.1

Let $X$ be the disjoint union of boundaries $\partial L_{H_v}$,
indexed by the vertices of $B$, and $p\colon X\to \vert B\subset B$
the obvious map. Then the boundary $\partial L$ is the tree completion
$\TC X$. The analogous statement holds also for rehabilitated
boundaries.  \endth

\pf Proof

We first describe $B$, following Serre \cite{Se}.  It has vertices
$\coprod_{V\in\vert Y}G/G_V$. We denote the vertex determined by $V$
and $[g]\in G/G_V$ by $g\widetilde V$. For each $E\in\edge Y$ and each
$[g]\in G/\init F_E$ there is an edge, denoted $g\widetilde E$, from
$g\widetilde{\init E}$ to $gt_E\widetilde{\term E}$. (Thus the reverse
of the edge $g\widetilde E$ is the edge determined by
$gt_E\widetilde{E^{-1}}$; this is slightly different notation from
\cite{Se}.)  Serre shows that $B$ is a tree and the quotient by the
obvious action of $G$ is $Y$.

We sketch the proof that $B$ is a tree.  We can choose a base vertex
for $B$ as the vertex $v_0=1\widetilde V_0$, where $V_0$ is the base
vertex for $\cal Y$.  A normal form representation $h=g_0t_{E_1}g_1
\ldots t_{E_{m}}g_m$ with $E_m=E$ for an element of $G$, as defined
early in section 3, determines a path in $B$ from the base vertex
$v_0$ of $B$ to the vertex $v=h\widetilde{\term E}$.  This path
consists of the sequence of edges $g_0\widetilde{E_1},\ldots,
g_{m-1}\widetilde{E_m}$.  Thus $B$ is a connected graph. Moreover,
since we can right-multiply $h$ by an element of $G_{\term E}$ without
changing $v$, we can assume $g_m\in\term F_E$. Then $h\in\cal G_E$ and
$h$ is determined up to the right-action of $F_E$.  Thus we see that
$v$ is actually determined by $E$ and an element of $\cal G_E/F_E$.
Now it not hard to see that any path without back-tracking from $v_0$
to $v$ gives a normal form representation for this $h$, and uniqueness
of normal forms up to the operations mentioned in section 3 leads to
uniqueness of such paths, showing that $B$ is a tree.

This also shows that the vertices of $B$ are in one-one correspondence
with\break $\coprod_{E\in\extY}\cal G_E/F_E$.  From this point of view the
stabilizer of the vertex $v=h\widetilde{\term E}$ corresponding to
$h\in \cal G_E/F_E$ is $H_v=hG_{\term E}h^{-1}$, and the induced
language $L_{H_v}$ is equivalent to the language $L_h$ on $G_{\term
E}$ (see Remark at end of section 3).

Now suppose $L$ is an asynchronous automatic structure on $G$. Since
the boundary $\partial L$ only depends on the equivalence class of
$L$, we may assume that $L$ is chosen as in Lemma 3.5.  We may also
assume it is prefix-closed.  Note that if $R$ is an $L$-rational
subset of $G$ and $R'$ is its ``prefix-closure'' (i.e., $R'=\eval N'$,
where $N'$ is the prefix closure of the set $N$ of $L$-words
evaluating into $R$), then $R'$ lies in a bounded neighborhood of $R$,
so $\partial R'=\partial R$. As bound one may take the diameter of a
finite state automaton for $N$.  Thus, in discussing the topology on
$\partial L$ we need only consider ``prefix-closed'' rational subsets
of $G$.

For any word $u\in L$, the shortest edge path decomposition (see
definition preceding Lemma 3.6) determines a shortest normal form
representative for $\eval u$, and hence, as above, a simple path
$\gamma_{\eval u}$ from the base vertex $v_0$ in $B$.  If $u_1$ is a
subword of $u$ then $\gamma_{\eval{u_1}}$ is a subpath of
$\gamma_{\eval u}$. It follows that an $L$-ray $w$ determines a simple
path $\gamma_w$ in $B$, which is a finite or infinite path according
as longer and longer initial segments of $w$ eventually all evaluate
into a fixed $hG_V$ or not.  We shall need the following Lemma.

\th Lemma 5.2

{\bf 1.}~~For any $k>0$ there exists $K>0$ such that if $u,u'\in
L$ satisfy $d(\eval u,\eval{u'})\le k$ then, by deleting at most the
last $K$ letters from $u$ and $u'$ one may obtain words $u_0$ and
$u'_0$ with $\gamma_{\eval{u_0}} = \gamma_{\eval{u'_0}}$.

{\bf 2.}~~If the $L$-ray $w$ fellow travels a subset $S\subset G$ then
every initial segment of $\gamma_w$ appears as an initial segment of
some $\gamma_g$, $g\in S$. The converse holds if $\gamma_w$ is
infinite and $S=\eval R$ with $R\subset L$ prefix-closed.\endth

\pf Proof

{\bf 1.}~ Since $L$ is finite-to-one, there exists a function
$\phi\colon\bbN\to\bbN$ such that any terminal segment $u_1$ of an
$L$-word with $\len(u_1)>\phi(k)$ satisfies $d(\eval{u_1},1)>k$.  Let
$k_1=\max\{d(f,1): f\in F_E$ for some edge $E$ of $\cal Y\}$ and
define $K=\phi(k+k_1)$.  Now suppose that $u,u'\in L$ satisfy $d(\eval
u,\eval{u'})\le k$ but do not satisfy the conclusion of the lemma.
Let $\gamma$ be the longest common segment of $\gamma_{\eval u}$ and
$\gamma_{\eval{u'}}$, and let $u_0$ and $u'_0$ be the longest initial
segments of $u$ and $u'$ with $\gamma_{\eval{u_0}}=
\gamma_{\eval{u'_0}}=\gamma$. Write $u=u_0u_1$, $u'=u'_0u'_1$.  At
least one of $u_1$ and $u'_1$, say $u_1$, has length greater than $K$.
By choice of $K$, the distance of $\eval{u_0}$ to $\eval
u=\eval{u_0u_1}$ exceeds $k+k_1$, so the distance from $\eval u$ to
$\eval{u_0}f$ exceeds $k$ for any $f$ in an edge group.  But, by
considering the normal form of $\eval{u^{-1}u'}$ one sees that the
shortest path in the Cayley graph from $\eval {u'}$ to $\eval u$ must
pass through $\eval{u_0}f$ for some $f\in\init F_E$, where $E$ is the
first edge of $u_1$ in the edge path decomposition of $u=u_0u_1$.
This path hence has length exceeding $k$, contradicting $d(\eval
u,\eval{u'})\le k$.

{\bf 2.}~ Part 1 of the lemma shows that if $w$ fellow travels a
subset $S$ of $G$ then there exist arbitrarily long initial segments
$u$ of $w$ with $\gamma_{\eval u}$ equal to an initial segment of a
path $\gamma_g$ with $g\in S$.  But every initial segment of
$\gamma_w$ is an initial segment of some such $\gamma_{\eval u}$, so
the first sentence of Lemma 5.2.2 is proved.

Conversely, suppose $\gamma_w$ is infinite and $S=\eval R$ with
$R\subset L$ prefix-closed and suppose every initial segment of
$\gamma_w$ is an initial segment of some $\gamma_g$ with $g\in S$.
For a given initial segment $\gamma$ of $\gamma_w$, choose such a
$g=\eval u$ with $u\in R$ and let $u_0$ and $w_0$ be the initial
segments of $u$ and $w$ corresponding to $\gamma$. Then $\eval {u_0}$
and $\eval{w_0}$ differ by an element of the edge group for the final
edge of $\gamma$, so $u_0$ and $w_0$ fellow travel.  Since $u_0$
travels in $S$ and $w_0$ is an arbitrarily long initial segment of
$w$, the result follows.\qed

We return to the proof of Theorem 5.1 for the boundary $\partial L$.
If two rays $w$ and $w'$ fellow travel then their paths $\gamma_w$ and
$\gamma_{{w'}}$ in $B$ are equal by Lemma 5.2.1.  Moreover, if
$\gamma_{w}$ and $\gamma_{{w'}}$ are equal and infinite then $w$ and
$w'$ do fellow travel by Lemma 5.2.2.  Thus rays $w$ with $\gamma_{w}$
infinite determine a subset of $\partial L$ that bijects to the set
$C$ of infinite rays in $B$.  This is the same as the set of ends of
$B$.

Suppose now $\gamma_{w}$ is finite, say it ends at the vertex of $B$
determined by $[h]\in\cal G_E/F_E$. Then by cutting $w$ at the point
where it has determined the whole path $\gamma_{w}$, we write $w$ in
the form $w_0u$ with $\eval{w_0}=h\eval f$ for some $f\in A_E$, and
$fu$ a ray in $L_h$.  If $w'$ is another ray with the same path
$\gamma_{w}$ we decompose it likewise as $w'_0u'$ with
$\eval{w'_0}=h\eval{f'}$ so that $f'u'$ is a ray in $L_h$.  Then $w$
and $w'$ fellow travel if and only if the rays $fu$ and $f'u'$ fellow
travel.  We thus get a copy of $\partial L_h$ in $\partial L$. We have
thus shown that, as a set, $\partial L$ is as claimed in the theorem.

It remains to verify that the topology is correct. Consider $v\in
\vert B$.  Let the simple path in $B$ from $v_0$ to $v$ be $\gamma$. 
The set $R_v:=\{u\in L:\gamma_{\eval u}=\gamma\}$ is regular. Hence,
$\eval{R_v}\subset G$ and its complement are both rational.  The rays
which fellow-travel $\eval{R_v}$ define the image of $\partial
L_{H_v}$ in $\partial L=X\cup C$ and the rays that fellow travel the
complement of $\eval{R_v}$ define the complement of $\partial
L_{H_v}$.  It follows that $\partial L_{H_v}$ is an open and closed
subset of $\partial L$. Moreover, rational subsets of $\eval{R_v}$
correspond to rational subsets of $H_v$, so $\partial L_{H_v}$ carries
the appropriate topology as a subspace of $\partial L$.  

Now suppose $w$ is a ray for which $\gamma_w$ is infinite, so $w$
represents an element of $C\subset \partial L$.  Suppose $S=\eval
R\subset G$ is a rational subset with $R\subset L$ prefix-closed, and
suppose $[w]$ is in the set $U$ of equivalence classes of rays that
fail to fellow-travel $S$.  Then Lemma 5.2.2 implies that there is
some initial segment $\gamma$ of $\gamma_w$ which does not appear as
an initial segment of any $\gamma_g$, $g\in S$.  Let $S_\gamma$ be the
set of $g\in G$ such that $\gamma_g$ has $\gamma$ as an initial
segment and let $U_\gamma\subset \partial L$ be the set of equivalence
classes of rays which fail to fellow-travel $G-S_\gamma$.  Then,
$[w]\in U_\gamma\subset U$, so these sets $U_\gamma$ form a
neighborhood basis for $[w]\in\partial L$. But $U_\gamma$ is the set
of equivalence classes of rays $w$ such that $\gamma_w$ has $\gamma$
as an initial segment. This defines the topology on $\partial L$
claimed in the theorem.

To see the analogous statement for the rehabilitated boundary we must
show that every non-trivial abstract simplex in $\partial L$ is an
abstract simplex of some $\partial L_{H_v}\subset \partial L$ and vice
versa. It is easy to see that an abstract simplex of $\partial
L_{H_v}$ is one for $\partial L$. Thus, we must show that if $x,y$ are
points of $\partial L$ which do not lie in some common $\partial
L_{H_v}$, then they do not form an abstract simplex, that is, they
have disjoint neighborhoods in $\cl G$.  By what was said above, a set
of the form $\cl \eval{R_v}$ is an open and closed subset of $\cl G$
whose intersection with $\partial L$ is $\partial L_{H_v}$. These sets
are disjoint for different $v$'s, so they provide disjoint
neighborhoods for $x$ and $y$ lying in distinct sets
$\partial{L_{H_v}}$.  Suppose just one of $x$ and $y$ lies in a
$\partial{L_{H_v}}$, say $x\in\partial{L_{H_v}}$ and $y\in C$. Then
$y=[w]$ with $\gamma_w$ infinite, so we can choose an initial segment
$\gamma$ of $\gamma_w$ which is not an initial segment of the path in
$B$ from $v_0$ to $v$.  The set $\cl S_\gamma$ is an open and closed
neighborhood of $y$ which is disjoint from the neighborhood
$\cl\eval{R_v}$ of $x$.  Finally, if $x=[w]$ and $y=[w']$ are distinct
points of $C$ and $\gamma$ and $\gamma'$ are initial segments of
$\gamma_w$ and $\gamma_{w'}$ which are longer than the longest common
initial segment of $\gamma_w$ and $\gamma_{w'}$ then $\cl S_\gamma$
and $\cl S_{\gamma'}$ are disjoint open and closed neighborhoods of
$x$ and $y$. \qed

\rk Remark

$\cl G$ can also be seen as a tree completion. For $g\in G$ the path
$\gamma_g$ ends in a vertex $v_g$ of $B$, so we get a map $p\colon
G\to\vert B\subset B$.  If $v$ is the vertex determined by $[h]\in\cal
G_{E}/F_{E}$, then $p^{-1}(v)=h(G_{\term E}-\term F_E)$ (except that
$p^{-1}(v_0)=G_{V_0}$ rather than $G_{V_0}-\{1\}$).  That is, up to a
finite set $p^{-1}(v)$ is just a translate of the group $G_{\pi(v)}$
on which the language $L_h$ (equivalent to $L_{H_v}$) is defined, so
we can attach to $p^{-1}(v)$ the boundary $\partial L_h\simeq\partial
L_{H_v}$.  This gives us a topology on $G\cup\coprod_{v\in\vert
B}\partial L_{H_v}$ and a map of this space to $\vert B\subset B$.
Then $\cl G$ is its tree completion. We leave the details to the
reader.

\Bib{References} 
 
\MaxReferenceTag{ECHLPT}
\rf{B} N. Brady, Asynchronous automatic structures on closed
hyperbolic surface groups, preprint, University of California 1993.

\rf{BGSS} G.  Baumslag, S. M.  Gersten, M.  Shapiro and H. Short, 
Automatic groups and amalgams,  Journal of Pure and
Applied Algebra {\bf 76} (1991), 229--316. 

\rf{ECHLPT}   D.B.A.  Epstein, J.W.  Cannon,   D.F.  Holt, S.V.F. Levy,
 M.S.  Paterson, and W.P.  Thurston, ``Word Processing in Groups,''
Jones and Bartlett Publishers, Boston, 1992.

\rf{GS} S. M.  Gersten and H.  Short, Rational sub-groups of 
biautomatic groups, Annals of Math. {\bf 134} (1991), 125--158. 

\rf{N} W. D. Neumann, Asynchronous combings of groups, Internat. J. 
Alg. Comp. {\bf 2} (1992), 179--185.

\rf{NS1} W. D. Neumann and M. Shapiro, Equivalent automatic structures
and their boundaries, Internat. J. Alg. Comp. {\bf 2} (1992),
443--469.

\rf{NS2} W. D. Neumann and M. Shapiro, Automatic structures and
geometrically finite hyperbolic groups, preprint.

\rf{Se} J. P. Serre, Trees, Graduate Texts in Mathematics (Springer
Verlag, 1980).

\rf{S1} M. Shapiro, Nondeterministic and deterministic asynchronous
automatic\break structures, Internat. J. Alg. and Comp. (1992). 

\rf{S2} M. Shapiro, Automatic structures and graphs of groups, in:
``Topology `90, Proceedings of the Research Semester in Low
Dimensional Topology at Ohio State,'' (Walter de Gruyter Verlag,
Berlin - New York 1992), 355--380.

\endBib

\Coordinates 
The Ohio State University\\ 
Department of Mathematics\\ 
Columbus, OH 43210 
\endCoordinates 

\Coordinates 
City College\\ 
Department of Mathematics\\ 
New York, NY 10031
\endCoordinates 

\bye